\documentclass[pdflatex,sn-mathphys-num]{sn-jnl}
\usepackage{graphicx}
\usepackage{multirow}
\usepackage{amsmath,amssymb,amsfonts}
\usepackage{amsthm}
\usepackage{mathrsfs}
\usepackage[title]{appendix}
\usepackage[dvipsnames]{xcolor}
\usepackage{textcomp}
\usepackage{manyfoot}
\usepackage{booktabs}
\usepackage{algorithm}
\usepackage{algorithmicx}
\usepackage{algpseudocode}
\usepackage{listings}
\usepackage{fullpage}
\usepackage{array}
\usepackage{Matharu_commands}
\usepackage{cleveref}
\crefname{equation}{}{}
\crefname{appendix}{Appendix}{Appendices}
\Crefname{appendix}{Appendix}{Appendices}
\crefname{theorem}{Remark}{Remarks}
\Crefname{theorem}{Remark}{Remarks}
\crefname{section}{Section}{Sections}
\Crefname{section}{Section}{Sections}
\crefname{subsection}{Section}{Sections}
\Crefname{subsection}{Section}{Sections}
\crefname{table}{Table}{Tables}
\Crefname{table}{Table}{Tables}
\usepackage{mathtools}
\usepackage{subcaption}

\algrenewcommand\algorithmicrequire{\textbf{Input:}}
\algrenewcommand\algorithmicensure{\textbf{Output:}}

\theoremstyle{thmstyletwo}

\newtheorem{remark}{Remark}

\theoremstyle{thmstyleone}
\newtheorem{lemma}{Lemma}%

\theoremstyle{thmstylethree}

\begin{document}

\title{Tolerance-driven close evaluation of the Stokes double layer potential on axisymmetric surfaces}

\author*[1]{\fnm{David} \sur{Krantz}}\email{davkra@kth.se}
\author[1,2]{\fnm{Pritpal} \sur{Matharu}}\email{matharu@mis.mpg.de}
\author[1]{\fnm{Anna-Karin} \sur{Tornberg}}\email{akto@kth.se}

\affil[1]{\orgdiv{Department of Mathematics}, \orgname{KTH Royal Institute of Technology}, \orgaddress{\street{Lindstedtsv{\"a}gen 25}, \city{Stockholm}, \postcode{11428}, \country{Sweden}}}
\affil[2]{\orgname{Max Planck Institute for Mathematics in the Sciences}, \orgaddress{\street{Inselstra{\ss}e 22}, \city{Leipzig}, \postcode{04103}, \country{Germany}}}

\abstract{We consider boundary integral methods for Stokes mobility and resistance problems involving smooth axisymmetric particles. A primary numerical challenge is the accurate and efficient evaluation of layer potentials at off-surface points close to particle surfaces. We present a tolerance-driven workflow for evaluating the Stokes double layer potential at such evaluation points (targets) to prescribed accuracy while avoiding unnecessary computational cost. For each target--particle interaction, a fast classifier selects the least costly option estimated to meet the tolerance among standard, upsampled, and special quadrature. Geometry-dependent unit-density error indicators are precomputed and tabulated in reduced cylindrical coordinates, then combined on the fly with a local layer-density modifier, making its cost negligible relative to evaluating the potential. Targets requiring special quadrature are treated using a stabilized version of singularity swap surface quadrature: the periodic azimuthal integral is evaluated first using translated singularity swap quadrature to prevent severe cancellation near the surface, followed by adaptive Gauss--Legendre quadrature in the polar direction guided by error indicators. We integrate this workflow into a boundary integral solver with precomputed quadrature by expansion for on-surface self-interactions and demonstrate the workflow’s performance for challenging configurations of spheroidal particles. Numerical results show that target classification is highly accurate. The prescribed tolerance is met for nearly all target--particle interactions and the error remains within a modest factor of the tolerance in the few remaining cases. Although the experiments focus on the Stokes double layer potential for spheroids, the off-surface framework applies to general smooth axisymmetric surfaces and can be easily adapted to other Stokes layer potentials.}





\maketitle

\section{Introduction}\label{s:introduction}

Particles suspended in viscous fluids arise throughout science and engineering, including microorganism locomotion \cite{Guasto_ARFM2012} and blood flow \cite{Zhao_JCP2010}. At low Reynolds numbers, the flow can be modeled by the Stokes equations. In this work, we consider Stokes flow exterior to a collection of smooth rigid axisymmetric particles.

Numerical approaches to particulate Stokes flow include, among others, immersed boundary methods \cite{Peskin_2002}, multiblob methods \cite{BROMS2023Multiblob}, and lattice-Boltzmann methods \cite{Wu2010Particle}.
For high-accuracy simulations, numerical methods based on boundary integral equations are particularly attractive. Since only the particle surfaces need to be discretized, no volume mesh is required, which is advantageous for moving particles and configurations with small interparticle gaps. Highly accurate boundary integral solvers have been developed for rigid-particle suspensions, including periodic and confined configurations \cite{AFKLINTEBERG2016,bagge2021highly,YAN2020109524,li2026scalableewaldfreebieframework}, as well as for slender bodies \cite{MALHOTRA2024112855}. A central numerical challenge is, however, the accurate evaluation of singular and nearly singular layer potentials, which generally requires special quadrature methods. For spheres and spheroids, this difficulty can be circumvented by expanding the layer density in spherical or spheroidal harmonics, respectively, for which the boundary integral operators and their off-surface potentials admit analytic representations \cite{CORONA2018327,CROWDER2026115040}. An alternative is the method of fundamental solutions (MFS), in which the flow is represented by singularities placed inside the particles, with their strengths determined by enforcing the boundary conditions. This avoids the need for singular and nearly singular surface quadrature. Recent MFS formulations treat mobility problems for large collections of ellipsoids \cite{Broms2025MFS}, while lubrication-adapted image sources enable accurate treatment of close interactions between spheres \cite{BROMS2025113636}.

We adopt a boundary integral equation formulation, reducing the Stokes problem to an integral equation posed on the particle surfaces \cite{pozrikidis1992boundary}. In the formulation considered here, the velocity is represented using the Stokes double layer potential, with the unknown layer density determined from the boundary integral equation. Our focus is the \emph{close evaluation} problem: the accurate and efficient evaluation of this potential at off-surface evaluation points (targets) located near particle surfaces. Such evaluations arise both in the solution of the integral equation, through interactions between nearby particles, and in the subsequent evaluation of the flow field close to the surfaces.

\begin{figure}[t]
\centering
\begin{subfigure}[t]{0.49\textwidth}
\includegraphics[width=\linewidth]{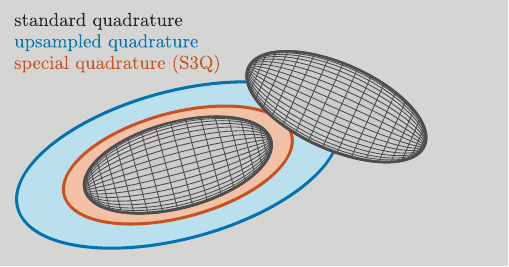}
\caption{}
\label{fig:quad_region_left}
\end{subfigure}
\hfill
\begin{subfigure}[t]{0.49\textwidth}
\includegraphics[width=\linewidth]{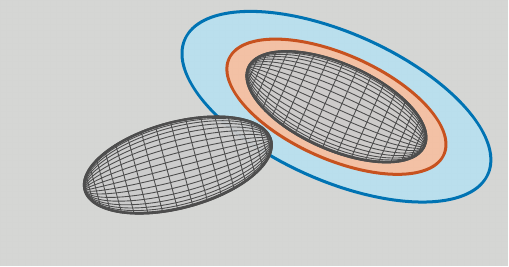}
\caption{}
\label{fig:quad_region_right}
\end{subfigure}
\caption{Illustration of quadrature regions around two axisymmetric particles. Panels (a) and (b) show the regions associated with the left and right spheroid, respectively. For fixed geometry, discretization, and layer density, the prescribed tolerance and target location determine which quadrature method should be used.}
\label{fig:quadrature_region_illustration}
\end{figure}

For an off-surface target, the integrand of the double layer potential is nonsingular but becomes sharply peaked as the target approaches a particle surface. The quadrature associated with the surface discretization, which we call standard quadrature, then loses accuracy. Assuming that the geometry and layer density are well resolved on this discretization, accuracy can initially be improved by interpolating the density to a finer grid and applying the same quadrature rule. We refer to this as upsampled quadrature. However, maintaining accuracy requires increasingly large upsampling factors as the target approaches the surface, eventually making this approach computationally impractical. At sufficiently small distances, a special quadrature method is therefore needed. Standard quadrature, upsampled quadrature with different factors, and special quadrature thus offer different levels of accuracy at different computational costs. For a prescribed tolerance, the appropriate strategy depends on the target location, particle geometry, and layer density and must be selected separately for each target--particle interaction. The resulting quadrature regions around two nearby spheroids are illustrated in \Cref{fig:quadrature_region_illustration}, where the individual upsampling levels are not distinguished.

Accurate quadrature for singular and nearly singular layer potentials in three dimensions remains an active area of research. Existing approaches include kernel regularization based on smooth kernel modifications \cite{Beale2026,faria2026highorderkernelregularizationsingular} or density interpolation \cite{FARIA2021113703}, local corrections to standard quadrature rules \cite{nitsche2025correctedtrapezoidalrulesnearsingular}, interpolation or extrapolation methods \cite{YING2006,MORSE2021,Bagge2023_thesis}, and methods that reduce surface integrals to line integrals \cite{zhu2022,jiang2024}. Related high-order solvers for axisymmetric surfaces use Fourier--Nyström schemes to exploit rotational symmetry, reducing the surface integral equation to modal equations on the generating curve and applying specialized quadrature to the resulting kernels \cite{YOUNG20124142,HELSING2014686}. Unlike these modal reductions, our approach requires only the individual particle surfaces to be axisymmetric; the overall particle configuration need not be rotationally symmetric.

Quadrature by expansion (QBX) is a unified approach for singular and nearly singular layer-potential evaluation~\cite{KLOCKNER2013,barnett2014}. It represents the potential through local expansions about centers placed away from the surface. For axisymmetric rigid particles, rotational symmetry enables efficient geometry-dependent precomputation, reducing on-surface evaluations at fixed nodes to the application of target-specific quadrature weights \cite{AFKLINTEBERG2016,bagge2021highly}. For arbitrary off-surface targets, the map from the layer density to the expansion coefficients can still be precomputed, but the application of the map and the target-dependent expansion evaluation must be performed online. We use QBX for on-surface evaluation, while the methods developed here address off-surface close evaluation.

For off-surface close evaluation, we build on singularity swap surface quadrature (S3Q) \cite{krantz2026s3q}, an extension of singularity swap quadrature (SSQ). SSQ was originally developed for nearly singular line integrals over open curves in two and three dimensions~\cite{AFKLINTEBERG2021} and subsequently extended to closed curves, first in two dimensions \cite{afKlinteberg2024,bao2024} and later in three dimensions \cite{krantz2026s3q}. In SSQ, the distance function in the denominator is paired with a simpler function having the same zeros in the complexified parameter plane. Multiplying and dividing by this function isolates the near singularity in the simpler factor, while the remaining ratio is smooth and can be combined with the layer density and expanded in an appropriate basis. The resulting basis integrals are then evaluated semi-analytically. For axisymmetric surfaces, S3Q applies SSQ to a trapezoidal-rule discretization in the azimuthal direction and evaluates the remaining polar integral on adaptively refined Gauss--Legendre panels, on which SSQ can also be applied. The required polar refinement and other target-dependent parameters are selected automatically from a user-prescribed error tolerance. Because S3Q relies on SSQ, it is, however, susceptible to catastrophic cancellation when the kernel numerator nearly vanishes, as occurs for the stresslet at targets very close to surface. Translated SSQ (TSSQ) addresses this instability using a target-dependent translated basis adapted to the local vanishing structure of the numerator \cite{krantz2025stabilizingsingularityswapquadrature}.

For large systems, i.e.~many particles, all interactions are initially evaluated using standard quadrature, with the resulting sums typically accelerated by fast summation methods, including fast multipole methods (FMMs) \cite{greengard1988,Greengard1987,Tornberg2008,Ying2004}, fast Ewald methods 
\cite{Lindbo2010,BAGGE2023}, and dual-space multilevel kernel splitting (DMK) \cite{jiang2025cpam,StokesDMK2026}. For each target--particle interaction for which standard quadrature is of insufficient accuracy, the corresponding contribution must then be corrected by locally replacing it with a more accurate evaluation. Reliable and efficient identification of these interactions is therefore essential, since missed corrections compromise accuracy, whereas unnecessary corrections add computational cost.

Previous Stokes particle solvers have selected among standard, upsampled, and special quadrature using distance-based regions, with thresholds calibrated numerically for a given geometry, discretization, and tolerance \cite{bagge2021highly}. 
Such criteria use the target-to-surface distance as a proxy, whereas the indicators derived in \cite{AFKLINTEBERG2022} and later refined for surfaces of spherical topology in \cite{SORGENTONE2023} directly approximate the quadrature error. The original S3Q workflow used these indicators to automatically select between standard quadrature and S3Q \cite{krantz2026s3q}. Since they can also approximate the error of upsampled quadrature, they provide a natural basis for extending this selection to several upsampling levels. Their direct evaluation, however, requires target-dependent computations, including complex root finding and the evaluation of geometric quantities, which can make the classification itself a non-negligible computational cost.

\subsection{Contributions and outline}\label{ss:contributions_outline}

Our aim is to provide a tolerance-driven and computationally efficient framework for evaluating Stokes double layer potentials at off-surface targets, including targets extremely close to a particle surface.
For each target--particle interaction, the framework uses standard quadrature when it is sufficiently accurate, otherwise selects the smallest prescribed upsampling factor predicted to meet the tolerance, and resorts to special quadrature only if no prescribed upsampling level suffices.
This requires error indicators that reliably distinguish between these regimes while adding only a negligible cost relative to the quadrature evaluation itself. 
To realize this framework, we make two main contributions.

Our first contribution is a fast procedure for evaluating the quadrature error indicators that drive this classification for smooth axisymmetric particles.
Building on the indicators developed in \cite{SORGENTONE2023}, we precompute and tabulate the required geometry-dependent quantities over a grid of target locations. Axisymmetry allows this grid to be represented using reduced target coordinates, substantially reducing the storage required for the resulting tables. During evaluation, interpolated values from these tables are combined with local information about the layer density, allowing the indicators to be evaluated at very low computational cost.

Our second contribution concerns the special quadrature component of the framework, which is based on S3Q \cite{krantz2026s3q}. Once special quadrature has been selected, S3Q automatically determines its internal parameters based on the prescribed tolerance, providing adaptive error control. Whereas previous demonstrations of S3Q used analytically prescribed layer densities, here we integrate it into a Stokes boundary integral solver and apply it both to close particle interactions when solving the integral equation and in subsequent post-processing. 
To make repeated S3Q evaluations more efficient, we precompute tables of SSQ quadrature weights, from which target-specific weights are obtained by interpolation.
To maintain the prescribed accuracy at extremely small target-to-surface distances, we incorporate TSSQ \cite{krantz2025stabilizingsingularityswapquadrature} into the azimuthal treatment and refer to the resulting variant as stabilized S3Q. This stabilization is particularly important for the Stokes double layer potential, since the near-vanishing of the stresslet numerator can otherwise cause severe numerical cancellation and loss of accuracy as the target approaches the surface.

Together, the fast target classifier and stabilized S3Q provide a tolerance-driven framework for off-surface evaluation on general smooth axisymmetric surfaces. We test the classifier on a capsule-shaped geometry, while the complete boundary-integral workflow, including precomputed QBX for on-surface self-interactions, is demonstrated for spheroids, for which this QBX implementation is available. Although the present work focuses on the Stokes double layer potential, the off-surface framework can be applied directly to other Stokes layer potentials, since the corresponding S3Q formulations and quadrature error indicators are already available.

\begin{figure}[t]
\centering
\includegraphics[width=0.85\linewidth]{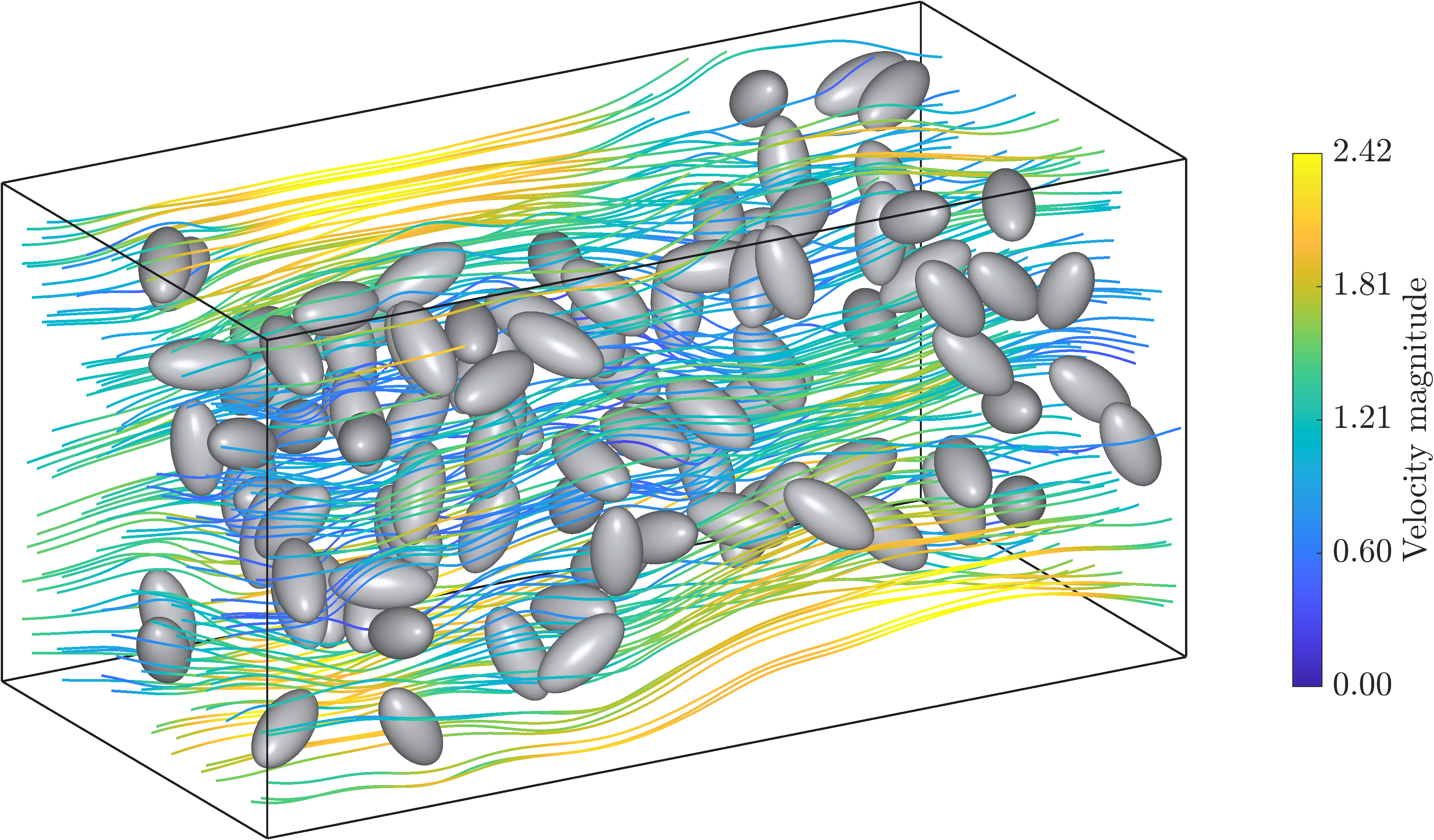}
\vspace{1em}
\caption{Triply periodic Stokes flow through 100 fixed spheroids. The streamlines are colored by the velocity magnitude.}
\label{fig:cluster_streamlines}
\end{figure}

The remainder of this paper is organized as follows. \Cref{s:problem_setting} introduces the Stokes boundary integral formulation, the axisymmetric surface discretization, standard and upsampled quadrature, and the QBX treatment of on-surface interactions. \Cref{s:precomputed_error_indicators} reviews the quadrature error indicators and develops their fast precomputed evaluation for target classification. \Cref{s:close_evaluation} presents the S3Q treatment, including the precomputation and interpolation of SSQ quadrature weights and the stabilization of the azimuthal evaluation using TSSQ. 
\Cref{s:numerical_experiments} assesses the accuracy and error control of the proposed framework through several numerical experiments. The final experiment applies the complete workflow to a model problem for porous flow in which 100 fixed spheroids are placed in a rectangular cell with periodic boundary conditions in all three coordinate directions; see \Cref{fig:cluster_streamlines}. In this example, the standard-quadrature contributions are accelerated using periodic DMK \cite{krantz2026dmk}.
Finally, \Cref{s:conclusions} presents the conclusions. Technical details from previous work are collected in Appendices \ref{a:error_formulas} and \ref{a:roots}.

\section{Problem setting}\label{s:problem_setting}

\subsection{The Stokes problem and boundary integral formulation}

We consider $m$ rigid axisymmetric particles with boundaries $\Gamma_{q} \subseteq \partial\Omega$, $q=1,\dots,m$, and write
\begin{equation*} \label{eq:boundary}
\Gamma \coloneqq \bigcup_{q=1}^{m} \Gamma_{q}
\end{equation*}
for the union of all particle surfaces. Let $\Omega\subset\RR^3$ denote the fluid domain exterior to the particles. On each $\Gamma_q$, the unit normal $\widehat{\vec{n}}$ is taken to point outward from the particle and into the fluid domain.
After nondimensionalization such that the viscosity is one, the velocity field $\u(\x)$ and pressure $p(\x)$ satisfy the Stokes equations
\begin{subequations}\label{eq:Stokes}
\begin{alignat}{2} 
\bfgrad p(\x) - \bfDelta \u(\x) &= \vec{0}, &  \qquad &\vec x\in\Omega,  \label{eq:StokesA} \\
\bfgrad \cdot \u(\x) &= 0, & \qquad &\vec x \in \Omega, \label{eq:StokesB}
\end{alignat} 
\end{subequations}
We allow for a prescribed background flow $\ubg$, which is a solution to \eqref{eq:Stokes} in the absence of the particles.
Unless stated otherwise, we consider the unbounded exterior problem and require the disturbance velocity $\u-\ubg$ to decay at infinity.
The no-slip boundary condition on each particle $q$ prescribes rigid-body velocity,
\begin{equation}\label{eq:no_slip_rigid_body}
\u(\x) = \u_q + \boldsymbol{\omega}_q \times (\x-\x_q),\qquad \x\in\Gamma_q,
\end{equation}
where $\u_q$, $\boldsymbol{\omega}_q$, and $\x_q$ denote the translational velocity, angular velocity, and center of particle $q$, respectively. 

We represent the disturbance flow using a completed Stokes double layer potential. With Einstein summation over repeated indices, the contribution from a surface $\Gamma_q$ is
\begin{align}\label{eq:dbl}
\D_i[\bsigma; \Gamma_q](\x) = \int_{\Gamma_q} \T_{ijk}(\x-\y)\sigma_j(\y) \widehat{n}_k (\y)\dif S(\y), \qquad i = 1, 2, 3, 
\end{align}
where $\bsigma$ is the unknown density, and $\bT$ is the stresslet kernel,
\begin{align}\label{eq:stresslet}
\T_{ijk}(\rr) = -6 \frac{r_i r_j r_k}{|\rr|^5}. 
\end{align}

The double layer potential alone cannot represent flows with nonzero net force or torque on the particles. To account for these, we add a so-called completion flow, introduced in \cite{power1987second}.
We place $N_{\mathrm{src}}$ completion sources at points $\y_{s,q}$ inside each particle $q$, giving
\begin{align} \label{eq:compflow}
\mathcal{V}_i(\x) = \sum_{q=1}^m \frac{1}{N_{\mathrm{src}}} \sum_{s=1}^{N_{\mathrm{src}}} \frac{1}{8\pi} \bigg(\mathcal{S}_{ij}(\x-\y_{s,q}) \,f_{j,q} + \mathcal{R}_{ij}(\x-\y_{s,q}) \,{t_{j,q}} \bigg),\quad \x\in\Omega,
\end{align}
where $\boldsymbol{f}_{q}$ and $\boldsymbol{t}_{q}$ are the net force and torque associated with particle $q$, respectively.
The Stokeslet and rotlet kernels are
\begin{align*} \label{eq:stokeslet_rotlet}
\mathcal{S}_{ij}(\rr) = \frac{\delta_{ij}}{|\rr|} + \frac{r_i r_j}{|\rr|^3},  \qquad \mathcal{R}_{ij}(\rr) = \varepsilon_{ijk} \frac{r_k}{|\rr|^3},
\end{align*}
with $\delta_{ij}$ the Kronecker delta symbol and $\varepsilon_{ijk}$ the Levi-Civita symbol. 

The velocity in the fluid is then represented as
\begin{equation}\label{eq:u_dlp}
\u(\x) = \ubg(\x) + \bD[\bsigma;\Gamma](\x) + \mathbfcal{V}(\x), \quad \x\in\Omega.
\end{equation}
Moreover, as $\x$ approaches the boundary, the double layer potential satisfies the jump relation
\begin{align} \label{eq:jump}
\underset{\varsigma \to 0}{\lim} \bD[\bsigma;\Gamma](\x\pm\varsigma\widehat{\n}) = \bD[\bsigma;\Gamma](\x) \mp 4\pi \bsigma(\x).
\end{align}
Taking the exterior limit in \eqref{eq:jump}, substituting it into the velocity representation, and imposing \eqref{eq:no_slip_rigid_body} yields the boundary integral equation
\begin{equation}
\u_q + \boldsymbol{\omega}_q \times (\x-\x_q)
= \ubg + \bD[\bsigma;\Gamma](\x) - 4\pi\bsigma(\x) + \mathbfcal{V}(\x), \qquad \x\in\Gamma.
\label{eq:integral_equation}
\end{equation}

We consider both mobility and resistance problems. In the mobility problem, the forces and torques $\boldsymbol{f}_{q}$ and $\boldsymbol{t}_{q}$ are prescribed, and hence the completion flow is known, while the translational and angular velocities $\u_q$ and $\boldsymbol{\omega}_q$ are unknown. These velocities can be expressed as functionals $\u_q[\bsigma]$ and $\boldsymbol{\omega}_q[\bsigma]$ of the layer density. In the resistance problem, the velocities are prescribed, while the forces and torques, and hence the completion flow, are instead expressed as functionals of $\bsigma$. Substituting the appropriate functionals into \eqref{eq:integral_equation} yields, in either case, a Fredholm boundary integral equation of the second kind for $\bsigma$; see, e.g., \cite[Section~2.1]{bagge2021highly} or \cite[p.~135]{pozrikidis1992boundary}.

We discretize \eqref{eq:integral_equation} using the Nyström method, enforcing the equation at the quadrature nodes on the particle surfaces. This yields a dense linear system for the nodal values of $\bsigma$, which is solved iteratively using the generalized minimal residual method (GMRES) \cite{saad1986}. Each application of the system matrix requires evaluating the double layer potential at the collocation points. 
Once $\bsigma$ has been computed, the unknown particle quantities---the rigid-body velocities in the mobility problem or the forces and torques in the resistance problem---are recovered from the corresponding functionals of $\bsigma$, and the velocity field in the fluid domain is evaluated using the representation formula \eqref{eq:u_dlp}.

\subsection{Surface discretization and quadrature}\label{ss:axisymmetric_discretization}
Working in a coordinate system centered on each particle and aligned with its symmetry axis, we parametrize the axisymmetric surface $\Gamma_q$ by the polar parameter $\conttheta\in[0,\pi]$ and azimuthal angle $\contphi\in[0,2\pi)$ as
\begin{equation}\label{eq:axisymmetric_param}
\boldsymbol{\gamma}(\conttheta,\contphi) = \big(a(\conttheta) \, \sin(\conttheta) \cos(\contphi),~a(\conttheta) \, \sin(\conttheta) \sin(\contphi),~c(\conttheta) \, \cos(\conttheta)\big),
\end{equation}
where $a(\conttheta)>0$ and $c(\conttheta)>0$ are smooth. For notational simplicity, we suppress the dependence of $\boldsymbol{\gamma}$, \(a\), and \(c\) on the particle index \(q\), as well as the rigid transformation that determines the particle's position and orientation in the relevant simulation setup.

The base surface discretization uses $n_\theta$ Gauss--Legendre nodes in the polar direction and $n_\phi$ equispaced trapezoidal nodes in the azimuthal direction. More generally, for an integer upsampling factor $\kappa\geq 1$, let $\{\theta_s^\kappa,w_s^{\theta,\kappa}\}_{s=1}^{\kappa n_\theta}$ denote the nodes and weights of the $\kappa n_\theta$-point Gauss--Legendre rule on $[0,\pi]$, and let $\{\phi_t^\kappa,w_t^{\phi,\kappa}\}_{t=1}^{\kappa n_\phi}$ denote the corresponding trapezoidal rule on $[0,2\pi)$.
To distinguish continuous surface parameters from quadrature nodes, we use $\conttheta$ and $\contphi$ for the continuous variables and $\theta$ and $\phi$ for the corresponding node coordinates.
We write surface functions with their pullbacks under $\boldsymbol{\gamma}$; e.g., $\bsigma(\conttheta,\contphi):=\bsigma\bigl(\boldsymbol{\gamma}(\conttheta,\contphi)\bigr)$.
For a target $\x\notin\Gamma_q$, the quadrature approximation of $\mathcal{D}_i[\bsigma;\Gamma_q]$ in \cref{eq:dbl} is then
\begin{align}
Q_i^{\kappa}[\bsigma;\Gamma_q](\boldsymbol{x})
=
\sum_{s=1}^{\kappa n_\theta}
\sum_{t=1}^{\kappa n_\phi}
&\mathcal{T}_{ijk}\!\left(
    \boldsymbol{x}
    -
    \boldsymbol{\gamma}
    \bigl(\theta_s^{\kappa},\phi_t^{\kappa}\bigr)
\right)
\sigma_j\bigl(\theta_s^{\kappa},\phi_t^{\kappa}\bigr)
\nonumber\\
&{}\times
\widehat{n}_k\bigl(\theta_s^{\kappa},\phi_t^{\kappa}\bigr)
\mathcal{J}\bigl(\theta_s^{\kappa},\phi_t^{\kappa}\bigr)
w_s^{\theta,\kappa}
w_t^{\phi,\kappa},
\qquad i=1,2,3.
\label{eq:upsampled_quadrature}
\end{align}
where $\mathcal{J}(\conttheta,\contphi)=\|\partial_\conttheta\boldsymbol{\gamma}\times\partial_\contphi\boldsymbol{\gamma}\|$ is the surface area element associated with the parametrization \cref{eq:axisymmetric_param}.
For $\kappa=1$, this is the quadrature associated with the original surface discretization, which we call standard quadrature. For $\kappa>1$, the density and the geometric grid functions appearing in \eqref{eq:upsampled_quadrature} are interpolated directly from the base grid to the finer $\kappa n_\theta\times\kappa n_\phi$ grid. We use barycentric Lagrange interpolation in the polar direction and trigonometric interpolation in the periodic azimuthal direction. We refer to the resulting approximation as upsampled quadrature with upsampling factor $\kappa$.

Upsampling assumes that the geometry and layer density are already sufficiently well resolved on the base grid. It introduces no additional degrees of freedom for $\bsigma$. Instead, the finer quadrature grid is used to resolve the increasingly rapid variation of the target-dependent kernel as the target approaches the surface. The number of quadrature nodes, and hence the direct evaluation cost, increases by a factor $\kappa^2$.

\subsection{QBX for on-surface evaluation}\label{ss:qbx}

QBX, introduced in \cite{KLOCKNER2013,barnett2014}, exploits the fact that layer potentials are smooth up to the boundary when approached from either side. They may therefore be represented by local expansions about centers placed off the boundary. Such an expansion can be evaluated at targets within its ball of convergence, including the boundary point where the ball touches the surface. For the Stokes double layer potential, we use the formulation of \cite{AFKLINTEBERG2016}, in which the potential is represented in terms of four Laplace dipole expansions.

Following \cite{AFKLINTEBERG2016,bagge2021highly}, for each surface node $\x_\ell\in\Gamma_q$, we place two expansion centers on the surface-normal line, one on the fluid side and one inside the particle. The corresponding expansions approximate the exterior and interior one-sided limits, denoted by $\vec\D^{+}$ and $\vec\D^{-}$, respectively. By the jump relation, their average equals the principal-value integral appearing in \cref{eq:integral_equation},
\begin{equation}
\label{eq:qbx_two_sided}
    \vec\D[\bsigma;\Gamma_q](\x_\ell)
    =
    \frac{1}{2}
    \left(
        \vec\D^{+}[\bsigma;\Gamma_q](\x_\ell)
        +
        \vec\D^{-}[\bsigma;\Gamma_q](\x_\ell)
    \right).
\end{equation}
We use this two-sided construction for all on-surface evaluations.

The main QBX parameters are the expansion radius, the truncation order, and the upsampling factor used to compute the expansion coefficients. 
The expansion radius determines the locations of the centers and the corresponding convergence balls, the truncation order determines the number of retained expansion terms, and the upsampling factor controls the quadrature accuracy of the coefficient integrals.
Together, these parameters control the achievable accuracy of QBX, and selecting them to meet a prescribed tolerance is nontrivial; see \cite[Section~6]{bagge2021highly} for a detailed discussion.

Let $N=n_\theta n_\phi$ denote the number of nodes on one particle. When solving the boundary integral equation, the QBX targets are the fixed nodes of the base surface discretization. For each target $\x_\ell$, the linear maps from the base-grid density values to the two one-sided QBX approximations can therefore be precomputed and combined according to \eqref{eq:qbx_two_sided} into a single set of target-specific quadrature weights. These weights incorporate the upsampling, computation of the expansion coefficients, and evaluation of the truncated expansions. Online evaluation then reduces to applying the weights to the current density, with the same $O(N)$ cost per target as direct quadrature. Since this cost is independent of the QBX parameters, they can be chosen conservatively to ensure high accuracy without increasing the cost of repeated evaluations, although this makes the precomputation more expensive.

Without exploiting symmetry, storing the target-specific weights for all $N$ nodes requires $O(N^2)$ memory, and applying them at all targets costs $O(N^2)$. Axisymmetry reduces the storage requirement because the weights for targets with the same polar coordinate are related by rotations and cyclic permutations of the surface data. It is therefore sufficient to store the weights for the $n_\theta$ targets on a single meridian, reducing the storage to $O(n_\theta N)$ \cite{AFKLINTEBERG2016}. If the particle is also reflection symmetric about its equatorial plane, only half of this meridian needs to be stored. 
The weights are precomputed once for each reference geometry and discretization for a given error level and can be reused for particles that differ only by rigid-body motion.

For arbitrary off-surface targets, the complete density-to-potential map cannot be precomputed because the target locations are not known in advance. Maps from the layer density to the expansion coefficients may still be precomputed, but these maps must be applied to the current density online, after which the resulting expansion must be evaluated at the target. Thus, upsampling only affects the precomputation, whereas the expansion order affects both storage and online evaluation cost.

In the present workflow, QBX is used only to evaluate the contribution from each particle at the nodes on its own surface, which we refer to as the self-contribution. Contributions from other particle surfaces, together with evaluations at targets in the fluid domain, are off-surface interactions and are handled by the quadrature-selection procedure developed in the next section.

\section{Fast precomputed error indicators and target classification}\label{s:precomputed_error_indicators}

In this section, we develop a target-classification strategy for off-surface target--particle interactions. For each target, smooth axisymmetric particle, and prescribed tolerance, the classifier determines whether standard quadrature is sufficiently accurate, one of the prescribed upsampling levels is needed, or special quadrature is required.
For this purpose, the error indicator need only approximate the magnitude of the quadrature error well enough to distinguish among these three regimes while remaining inexpensive to evaluate.

We begin by recalling the quadrature error analysis underlying the indicator. Sections \ref{ss:error_formulas_one_dimension} and \ref{ss:error_formulas_surfaces} consider scalar-valued potentials with a general singularity exponent $p$, first on curves and then
on surfaces. In Section \ref{ss:tensor_indicator}, we specialize the analysis to the tensor-valued stresslet kernel of the Stokes double layer potential, for which $p=5/2$. \Cref{ss:tabulation} describes the tabulation and fast evaluation of the resulting indicator, and \Cref{ss:performance_precomputed_indicators} assesses its accuracy and computational cost.

\subsection{Error indicators for line potentials on curves}\label{ss:error_formulas_one_dimension}

We begin with a scalar-valued potential obtained by integrating along a smooth curve, which we refer to as a line potential.
Let $\interval$ denote the reference base parameter interval, taken to be $[-1,1]$ for Gauss--Legendre quadrature and $[0,2\pi)$ for the trapezoidal rule, and $\vec y:\interval\subset\mathbb{R}\rightarrow\mathbb{R}^3$ be an analytic parametrization of the curve. For the trapezoidal rule, the parametrization and numerator are assumed
to be $2\pi$-periodic.

For a target point $\xtil=(x_1,x_2,x_3)\in\mathbb{R}^3$ lying off the curve, we define the squared-distance function associated with the curve $\vec y=(y_1,y_2,y_3)$ as
\begin{equation}
R_{\vec{y}}^2(t, \xtil) \coloneqq \sum_{i=1}^3 \left(y_i(t) - x_i\right)^2.
\label{eq:R2_curve}
\end{equation}

For a smooth numerator $g$, we define the integrand of a generic line potential
\begin{equation}
v_p(t) = \frac{g(t)}{\big(R_{\vec{y}}^2(t,\xtil)\big)^p}, \qquad  2p\in\mathbb{Z}^+.
\label{eq:vp_scalar}
\end{equation}
Let
\begin{equation}
I[v_p](\xtil) = \int_\interval v_p(t)\dif t, \qquad Q_n[v_p] = \sum_{i=1}^nv_p(t_i)w_i,
\label{eq:layer_pot_curve}
\end{equation}
denote its exact integral and $n$-point quadrature approximation, respectively, where $\{t_i\}_{i=1}^n$ and $\{w_i\}_{i=1}^n$ are the quadrature nodes and weights. The associated quadrature error is
\begin{equation}
E_{n}[v_p]\coloneqq \Big|I[v_p]-Q_n[v_p]\Big|.
\label{eq:En}
\end{equation}
For integrands of the form \eqref{eq:vp_scalar}, we also write $E_{n,p}[g] \coloneqq E_n[v_p]$,
to indicate the dependence on the numerator $g$ and the exponent $p$.
The dependence on the curve and target point is left implicit.
When needed, the quadrature rule is indicated by a superscript,
as in $E_{n,p}^{\GL}[g]$ and $E_{n,p}^{\TZ}[g]$.

The difference $I[v_p]-Q_n[v_p]$ can be represented as a contour integral in the complex plane involving the analytically continued integrand $v_p$ and a remainder function specific to the quadrature rule \cite{DONALDSON1972,ELLIOT2008}. The contour encloses the integration interval and quadrature nodes and lies within a region where $v_p$ is analytic. We assume that $g$ and $\vec y$ admit analytic continuations to a complex neighborhood of $\interval$ containing the relevant roots of $R_{\vec y}^2$. Within this neighborhood, singularities of $v_p$ can therefore occur only at these roots. Deforming the contour toward these singularities shows that those closest to the integration interval provide the leading contribution to the quadrature error.

We assume that the roots considered below are simple and that $g$
does not vanish there. The corresponding singularities of $v_p$
are poles of order $p$ when $p\in\ZZ^+$, and branch points otherwise.
Since $R_{\vec{y}}^2(t, \xtil)$ is real for $t\in \mathbb{R}$, its roots come in complex-conjugate pairs. Let
\begin{equation*}
t_0=t_0(\xtil),\quad t_0^*=\big(t_0(\xtil)\big)^*
\label{eq:t_root}
\end{equation*}
denote the complex-conjugate root pair closest to the real axis, satisfying
\begin{align*}
R_{\vec{y}}^2(t_0, \xtil) = R_{\vec{y}}^2(t_0^*, \xtil) = 0.
\end{align*}
In the following, we assume that this pair gives the dominant contribution to the quadrature error.

To separate out the factor that depends on the quadrature rule, we define
\begin{align}
\Psi_{n,p}^{\mathrm{TZ}}(t_0)
&:=
\frac{4\pi n^{p-1}}{\Gamma(p)}
e^{-n|\Im(t_0)|},
\label{eq:tz_root_factor}
\\
\Psi_{n,p}^{\mathrm{GL}}(t_0)
&:=
\frac{4\pi}{\Gamma(p)}
\left|
    \frac{2n+1}{\sqrt{t_0^2-1}}
\right|^{p-1}
\frac{1}{\varrho(t_0)^{2n+1}},
\label{eq:gl_root_factor}
\end{align}
where $\Gamma(p)$ is the gamma function, and $\varrho(t)$ denotes the Bernstein radius $\varrho(t)\coloneqq\big|t+\sqrt{t^2-1}\big|$,
where $\sqrt{t^2-1}$ is defined as $\sqrt{t+1}\sqrt{t-1}$, with branch cuts chosen so that $-\pi<\arg(t\pm1)\leq\pi$. 

Then, as shown in \cite[Section 5]{AFKLINTEBERG2022}, the quadrature error \eqref{eq:En} for the $n$-point trapezoidal rule and Gauss--Legendre rules can be approximated as
\begin{align}
E_{n,p}^{\TZ}[g]
&\approx
|g(t_0)|
|G(t_0)|^p
\Psi_{n,p}^{\TZ}(t_0),
\label{eq:etz_1d}
\\
E_{n,p}^{\GL}[g]
&\approx
|g(t_0)|
|G(t_0)|^p
\Psi_{n,p}^{\GL}(t_0),
\label{eq:egl_1d}
\end{align}
Here, for a simple root $w$ of $R^2_{\vec y}(\cdot,\xtil)$, the geometry factor is defined by
\begin{align}
G(w):=\lim_{t \rightarrow w} \frac{t-w}{R_{\vec{y}}^2(t, \xtil)}=\big( 
2 (\vec y(w)-\xtil)\cdot \vec y'(w)
\big)^{-1}.
\label{eq:geom_fac}
\end{align}
We do not explicitly indicate the dependence on $\xtil$, but $t_0$, $G$, and, depending on the kernel, possibly also $g$, all depend on it. The formulas above are asymptotic in the quadrature order $n$, and are therefore formally justified as $n\rightarrow\infty$. In practice, however, the resulting indicators are already useful for moderate values of $n$, and have been found to predict the loss of accuracy remarkably well \cite{AFKLINTEBERG2022}.

\begin{remark}
As noted above, $\{t_0,t_0^*\}$ is the singularity pair closest to the real axis. There may be additional singularity pairs, which also contribute to the quadrature error.
However, as seen from \cref{eq:tz_root_factor,eq:gl_root_factor}, the exponential
decay of each contribution is governed by $|\Im(t)|$ for the
trapezoidal rule and by $\varrho(t)$ for Gauss--Legendre quadrature.
For comparable numerator and geometry prefactors, pairs with larger
values of the corresponding measure therefore contribute
exponentially less as $n$ increases.
In practice, even when two singularity pairs make comparable
contributions to the indicator, retaining only the closest one
changes its value by only a modest constant factor, and can
therefore still capture the error magnitude sufficiently well
for target classification.
\end{remark}

\subsection{Error indicators for layer potentials on surfaces}\label{ss:error_formulas_surfaces}
We now consider a scalar-valued layer potential on the axisymmetric surface parametrized by $\bgamma(\conttheta, \contphi)$ in Section \ref{ss:axisymmetric_discretization}. For an off-surface target $\xtil$, define the squared-distance function
\begin{equation}
R_{\bgamma}^2(\conttheta, \contphi, \xtil) \coloneqq \sum_{i=1}^3 \left(\gamma_i(\conttheta, \contphi) - x_i\right)^2.
\label{eq:R2}
\end{equation}
A generic scalar-valued layer potential then takes the form
\begin{align}
U_p[f\sigma](\xtil) = \int_{0}^{\pi}\int_{0}^{2\pi} \frac{f(\conttheta, \contphi)\,\sigma(\conttheta,\contphi)}{\big(R_{\bgamma}^2(\conttheta, \contphi, \xtil)\big)^p}\dif\contphi\dif\conttheta,\qquad 2p\in\mathbb{Z}^+,
\label{eq:layerpotparam}
\end{align}
where $f$ contains the smooth kernel and parametrization factors, including the surface Jacobian, and $\sigma$ is the scalar layer density. Any target dependence of $f$ is left implicit.

The surface is discretized as introduced in Section \ref{ss:axisymmetric_discretization},
with $n_\theta$ Gauss--Legendre nodes in the polar direction and $n_\phi$ trapezoidal nodes in the azimuthal direction. Let $I_\conttheta$ and $I_\contphi$ denote integration in the respective variables, and let $Q_{n_\theta}^\GL$ and $Q_{n_\phi}^\TZ$ denote the corresponding quadrature rules. The polar rule is understood to include the mapping of nodes and weights to $[0,\pi]$.
Temporarily introducing $v_p$ as the integrand of \eqref{eq:layerpotparam}, the quadrature error of \eqref{eq:layerpotparam} is
\begin{equation*}
{E}_{n_\theta,n_\phi}[f\sigma]\coloneqq \left| \left( {I}_\conttheta {I}_\contphi - {Q}^{\GL}_{n_\theta} {Q}^{\TZ}_{n_\phi} \right)[v_p] \right| 
\approx \left({I}_\conttheta {E}^{\TZ}_{n_{\phi},p} + 
{I}_\contphi {E}^{\GL}_{n_{\theta},p}\right)[f\sigma],
\label{eq:2DerrIntro}
\end{equation*}
where we, following \cite{AFKLINTEBERG2022}, have neglected the quadratic error term. 
Hence, we can compute an approximation of the tensor-product quadrature error by integrating the one-dimensional error estimates that we introduced in the previous section.  
Specifically, define
\begin{align}
E_{n_\phi,p}^{\contphi}[f\sigma](\xtil)
&:=
\int_0^\pi
E_{n_\phi,p}^{\TZ}
\bigl[f(\conttheta,\cdot)\,\sigma(\conttheta,\cdot)\bigr]
\dif\conttheta,
\label{eq:surface_error_phi}
\\
E_{n_\theta,p}^{\conttheta}[f\sigma](\xtil)
&:=
\int_0^{2\pi}
E_{n_\theta,p}^{\GL}
\bigl[f(\cdot,\contphi)\,\sigma(\cdot,\contphi)\bigr]
\dif\contphi.
\label{eq:surface_error_theta}
\end{align}
Then their sum provides an indicator for the surface quadrature error,
\begin{equation}
E_{n_\theta,n_\phi}[f\sigma](\xtil)
\approx
E_{n_\phi,p}^{\contphi}[f\sigma](\xtil)
+
E_{n_\theta,p}^{\conttheta}[f\sigma](\xtil),
\label{eq:surface_error_sum}
\end{equation}
where each one-dimensional error in \cref{eq:surface_error_phi,eq:surface_error_theta} can be treated using the formulas from Section~\ref{ss:error_formulas_one_dimension}.

To apply the one-dimensional formulas, we consider the curves
obtained by holding one surface parameter fixed.
For fixed $\conttheta$, the azimuthal curve is
$\vec y(\contphi)=\bgamma(\conttheta,\contphi)$, while for fixed
$\contphi$, the polar curve is
$\vec y(\conttheta)=\bgamma(\conttheta,\contphi)$.
The squared-distance function \eqref{eq:R2_curve} along each curve
is obtained by fixing the corresponding parameter in
$R_{\bgamma}^2$.

Let $\contphi_0(\conttheta)$ and $\conttheta_0(\contphi)$ denote
the selected roots in the azimuthal and polar directions,
respectively. Together with their complex conjugates, they satisfy
\begin{align*}
R_{\bgamma}^2
\bigl(\conttheta,\contphi_0(\conttheta),\xtil\bigr)
&=
R_{\bgamma}^2
\bigl(\conttheta,\contphi_0^*(\conttheta),\xtil\bigr)
=0,
\\
R_{\bgamma}^2
\bigl(\conttheta_0(\contphi),\contphi,\xtil\bigr)
&=
R_{\bgamma}^2
\bigl(\conttheta_0^*(\contphi),\contphi,\xtil\bigr)
=0.
\end{align*}
In each direction, we select the pair closest to the real axis
and make the same simple-root and dominant-pair assumptions as in
Section~\ref{ss:error_formulas_one_dimension}.
The dependence of the roots on $\xtil$ is left implicit.

The geometry factor \eqref{eq:geom_fac} applies along each coordinate curve. To distinguish the
two directions, we define
\begin{align*}
G_\conttheta(\conttheta,\contphi)
&:=
\left[
    2\bigl(\bgamma(\conttheta,\contphi)-\xtil\bigr)
    \cdot
    \partial_{\conttheta}\bgamma(\conttheta,\contphi)
\right]^{-1},
\\
G_\contphi(\conttheta,\contphi)
&:=
\left[
    2\bigl(\bgamma(\conttheta,\contphi)-\xtil\bigr)
    \cdot
    \partial_{\contphi}\bgamma(\conttheta,\contphi)
\right]^{-1},
\end{align*}
where the subscripts indicate the direction of differentiation.

Applying \cref{eq:etz_1d,eq:egl_1d} gives the two directional contributions \eqref{eq:surface_error_sum} as
\begin{align*}
E_{n_\phi,p}^{\varphi}[f\sigma](\xtil)
&\approx
\int_0^\pi
\bigl|
    f\bigl(\vartheta,\varphi_0(\vartheta)\bigr)
    \sigma\bigl(\vartheta,\varphi_0(\vartheta)\bigr)
\bigr|
\bigl|
    G_\varphi\bigl(\vartheta,\varphi_0(\vartheta)\bigr)
\bigr|^p
\Psi_{n_\phi,p}^{\TZ}
\bigl(\varphi_0(\vartheta)\bigr)
\dif\vartheta,
\\
E_{n_\theta,p}^{\vartheta}[f\sigma](\xtil)
&\approx
\left(\frac{\pi}{2}\right)^{1-p}
\int_0^{2\pi}
\bigl|
    f\bigl(\vartheta_0(\varphi),\varphi\bigr)
    \sigma\bigl(\vartheta_0(\varphi),\varphi\bigr)
\bigr|
\bigl|
    G_\vartheta\bigl(\vartheta_0(\varphi),\varphi\bigr)
\bigr|^p
\Psi_{n_\theta,p}^{\GL}
\bigl(t_0(\varphi)\bigr)
\dif\varphi,
\end{align*}
with $\Psi_{n_\phi,p}^{\TZ}$ and $\Psi_{n_\theta,p}^{\GL}$ defined in \cref{eq:tz_root_factor,eq:gl_root_factor}. Here $\Psi_{n_\theta,p}^{\GL}$ is evaluated at the mapped polar
root $t_0(\contphi)=2\conttheta_0(\contphi)/\pi-1$, corresponding
to the reference interval $[-1,1]$. The factor $(\pi/2)^{1-p}$
accounts for the same change of variables.

The most rapid variation in these integrals typically comes
through the root-dependent factors $\Psi$.
To obtain an inexpensive indicator of the error magnitude, we
approximate the remaining factors by their values at reference
parameters near the target.
Let $(\bar\vartheta,\bar\varphi)$ denote the parameter values of
the surface point closest to $\xtil$\footnote{In practice, these reference parameters are approximated
by those of the closest surface grid point.}, and define the
reference roots
\begin{equation*}
    \bar\varphi_0
    :=
    \varphi_0(\bar\vartheta),
    \qquad
    \bar\vartheta_0
    :=
    \vartheta_0(\bar\varphi).
    \label{eq:surface_reference_roots}
\end{equation*}
Evaluating the numerator and geometry prefactors at these
reference values leaves only $\Psi$ under each integral.

For unit density, the directional indicators then become
\begin{align}
E_{n_\phi,p}^{\varphi}[f](\xtil)
&\approx
\bigl|f(\bar\vartheta,\bar\varphi_0)\bigr|
\bigl|G_\varphi(\bar\vartheta,\bar\varphi_0)\bigr|^p
\int_0^\pi
\Psi_{n_\phi,p}^{\TZ}
\bigl(\varphi_0(\vartheta)\bigr)
\dif\vartheta,
\label{eq:surface_unit_indicator_phi}
\\
E_{n_\theta,p}^{\vartheta}[f](\xtil)
&\approx
\left(\frac{\pi}{2}\right)^{1-p}
\bigl|f(\bar\vartheta_0,\bar\varphi)\bigr|
\bigl|G_\vartheta(\bar\vartheta_0,\bar\varphi)\bigr|^p
\int_0^{2\pi}
\Psi_{n_\theta,p}^{\GL}
\bigl(t_0(\varphi)\bigr)
\dif\varphi.
\label{eq:surface_unit_indicator_theta}
\end{align}
For a general density, the same approximation gives
\begin{align}
E_{n_\phi,p}^{\varphi}[f\sigma](\xtil)
&\approx
\bigl|\sigma(\bar\vartheta,\bar\varphi_0)\bigr|
E_{n_\phi,p}^{\varphi}[f](\xtil),
\label{eq:surface_density_scaling_phi}
\\
E_{n_\theta,p}^{\vartheta}[f\sigma](\xtil)
&\approx
\bigl|\sigma(\bar\vartheta_0,\bar\varphi)\bigr|
E_{n_\theta,p}^{\vartheta}[f](\xtil).
\label{eq:surface_density_scaling_theta}
\end{align}
This separation is useful because the unit-density indicators
depend on the kernel, geometry, discretization, and target point,
but not on the layer density. They can therefore be precomputed
and combined with local density information during evaluation.

Define the (local) \emph{density modifier} by
\begin{equation}
\sigma^M(\xtil)
:=
\max\left\{
    \bigl|\sigma(\bar\vartheta,\bar\varphi_0)\bigr|,
    \bigl|\sigma(\bar\vartheta_0,\bar\varphi)\bigr|
\right\}.
\label{eq:scalar_density_modifier}
\end{equation}
Combining \eqref{eq:surface_error_sum} with \cref{eq:surface_density_scaling_phi,eq:surface_density_scaling_theta,eq:scalar_density_modifier} then gives
\begin{equation*}
\begin{aligned}
    E_{n_\theta,n_\phi}[f\sigma](\xtil)
    \approx{}&
    \bigl|\sigma(\bar\vartheta,\bar\varphi_0)\bigr|
    E_{n_\phi,p}^{\varphi}[f](\xtil)+
    \bigl|\sigma(\bar\vartheta_0,\bar\varphi)\bigr|
    E_{n_\theta,p}^{\vartheta}[f](\xtil)
    \\
    \le{}&
    \sigma^M(\xtil)
    \left(
        E_{n_\phi,p}^{\varphi}[f](\xtil)
        +
        E_{n_\theta,p}^{\vartheta}[f](\xtil)
    \right).
\end{aligned}
\label{eq:surface_density_modified_indicator}
\end{equation*}

\begin{remark}[Practical evaluation of the indicator]\label{rem:practical_evaluation_indicator}
The azimuthal root $\contphi_0(\conttheta,\xtil)$ is available in closed form for any smooth axisymmetric surface, whereas the polar root $\conttheta_0(\contphi,\xtil)$ admits a closed-form expression only for special geometries such as the sphere and the spheroid; see Appendix~\ref{a:roots}. For a general axisymmetric surface, $\conttheta_0$ must therefore be computed by a one-dimensional root-finding procedure. 
To limit the number of root evaluations when approximating the integrals in \cref{eq:surface_unit_indicator_phi,eq:surface_unit_indicator_theta}, we follow \cite{AFKLINTEBERG2022} and use local approximations to the root functions, corrected using the reference values $\bar{\conttheta}_0$ and $\bar{\contphi}_0$. Consequently, only these two reference roots need to be computed explicitly for each target point. The resulting integrals are then evaluated using the low-order Gauss--Laguerre quadrature described in Appendix~\ref{a:error_formulas}.
\end{remark}

\subsection{Error indicators for the tensor-valued stresslet kernel}\label{ss:tensor_indicator}
We now explain how the scalar-valued layer potential error indicator construction of Section \ref{ss:error_formulas_surfaces} is used for the vector-valued Stokes double layer potential \eqref{eq:dbl}. Here we write sums over vector components explicitly, rather than using Einstein summation, in order to distinguish the output velocity component $i$ from the input density component $j$.

We begin by rewriting \eqref{eq:dbl} using the surface parametrization in Section \ref{ss:axisymmetric_discretization} for one particle surface $\Gamma_q$. By contracting the stresslet kernel \eqref{eq:stresslet} with the normal vector we find
\begin{equation}
\D_i[\bsigma; \Gamma_q](\xtil) = \sum_{j=1}^3  \int_{0}^{\pi}\int_{0}^{2\pi} \frac{f_{ij}(\conttheta,\contphi,\xtil)\,\sigma_j(\conttheta,\contphi)}{\big(R^2_{\bgamma}(\conttheta, \contphi,\x)\big)^{5/2}}\dif\contphi\dif\conttheta,\quad i=1,2,3,
\label{eq:dbl_sum}
\end{equation}
where $R_{\bgamma}^2$ is defined in \eqref{eq:R2}, and
\begin{equation}
f_{ij}(\conttheta,\contphi,\xtil) = -6r_ir_j(\vec r\cdot\widehat{\n})\mathcal{J}(\conttheta,\contphi),\qquad \vec r=\xtil-\bgamma(\conttheta,\contphi).
\label{eq:fij_r}
\end{equation}
Here $\widehat{\n}$ is the outward unit normal at $\bgamma(\conttheta,\contphi)$, and $\mathcal{J}$ is the surface Jacobian defined below \eqref{eq:upsampled_quadrature}.
For each fixed pair $(i,j)$, the integral in \eqref{eq:dbl_sum} has the scalar form \eqref{eq:layerpotparam}, with $p=5/2$, smooth factor $f_{ij}$, and scalar density $\sigma_j$.

Let $\vec e_j$, $j=1,2,3$, denote the constant unit density with $(\vec e_j)_\ell=\delta_{j\ell}$. If $\vec{\sigma}=\vec e_j$, then only the $j$th density component is nonzero, and \eqref{eq:dbl_sum} reduces to the scalar unit-density integral $\mathcal{D}_i[\vec e_j;\Gamma_q](\xtil)=U_{5/2}[f_{ij}](\xtil)$ with $U_{5/2}$ from \eqref{eq:layerpotparam}. We therefore define the directional unit-density indicators by
\begin{align}
    E_{ij}^{\contphi}(\xtil)
    &:=
    E_{n_\phi,5/2}^{\contphi}[f_{ij}](\xtil),
    \label{eq:stresslet_unit_indicator_phi}
    \\
    E_{ij}^{\conttheta}(\xtil)
    &:=
    E_{n_\theta,5/2}^{\conttheta}[f_{ij}](\xtil),
    \label{eq:stresslet_unit_indicator_theta}
\end{align}
and their sum by
\begin{equation*}
    E_{ij}(\xtil)
    :=
    E_{ij}^{\contphi}(\xtil)
    +
    E_{ij}^{\conttheta}(\xtil).
    \label{eq:stresslet_unit_indicator}
\end{equation*}
These quantities are evaluated using
\cref{eq:surface_unit_indicator_phi,eq:surface_unit_indicator_theta},
with $f$ replaced by $f_{ij}$.

For a general vector-valued density
$\vec{\sigma}=(\sigma_1,\sigma_2,\sigma_3)$, we apply the density
scaling in
\cref{eq:surface_density_scaling_phi,eq:surface_density_scaling_theta}
to each input component.
Analogously to \eqref{eq:scalar_density_modifier}, define
\begin{equation}
    \sigma_j^M(\xtil)
    :=
    \max\left\{
        \bigl|
            \sigma_j(\bar{\conttheta},\bar{\contphi}_0)
        \bigr|,
        \bigl|
            \sigma_j(\bar{\conttheta}_0,\bar{\contphi})
        \bigr|
    \right\},
    \qquad j=1,2,3.
    \label{eq:stresslet_density_modifier}
\end{equation}
Summing the contributions from the three input components gives
the output-component indicator
\begin{equation*}
    \begin{aligned}
        E_i[\vec{\sigma}](\xtil)
        :={}&
        \sum_{j=1}^3
        \Bigl(
            E_{ij}^{\contphi}(\xtil)
            \bigl|
                \sigma_j(\bar{\conttheta},\bar{\contphi}_0)
            \bigr|
            +
            E_{ij}^{\conttheta}(\xtil)
            \bigl|
                \sigma_j(\bar{\conttheta}_0,\bar{\contphi})
            \bigr|
        \Bigr)
        \\
        \le{}&
        \sum_{j=1}^3
        E_{ij}(\xtil)\sigma_j^M(\xtil),
        \qquad i=1,2,3.
    \end{aligned}
    \label{eq:stresslet_component_indicator}
\end{equation*}
Taking the maximum over the output components gives
\begin{equation*}
    E[\vec{\sigma}](\xtil)
    :=
    \max_{i=1,2,3}
    \sum_{j=1}^3
    E_{ij}(\xtil)\sigma_j^M(\xtil).
    \label{eq:stresslet_max_indicator}
\end{equation*}

To avoid storing and interpolating the full $3\times 3$ collection of unit-density indicators $E_{ij}$ at runtime, we collapse the output index, and define
\begin{equation}
    E_j(\xtil)
    :=
    \max_{i=1,2,3}E_{ij}^{\contphi}(\xtil)
    +
    \max_{i=1,2,3}E_{ij}^{\conttheta}(\xtil),
    \qquad j=1,2,3.
    \label{eq:stresslet_collapsed_indicator}
\end{equation}
Since $E_{ij}(\xtil)\le E_j(\xtil)$ and the modifiers
$\sigma_j^M$ are nonnegative, it follows that
\begin{equation}
    E[\vec{\sigma}](\xtil)
    \le
    \sum_{j=1}^3 E_j(\xtil)\sigma_j^M(\xtil).
    \label{eq:stresslet_classification_indicator}
\end{equation}
The right-hand side of \eqref{eq:stresslet_classification_indicator} is the error indicator used for target classification and requires only one unit-density indicator per input density component.

\subsection{Tabulation in reduced target coordinates}\label{ss:tabulation}
To speed up evaluation of the error indicator, we exploit axisymmetry to tabulate the unit-density indicators in reduced target coordinates. 
When the polar root is not available in closed form, we also tabulate it and recover it by interpolation. This additional root table is not needed for spheroids, but allows the same classifier to be applied to more general smooth axisymmetric geometries.
These tables are precomputed for a fixed particle geometry and discretization and can be reused for different targets and densities.

We work in particle-centered coordinates with the symmetry
axis along the third coordinate. For a target
$\xtil=(x_1,x_2,x_3)$, define
\begin{equation*}
\rho = \sqrt{x_1^2+x_2^2},
\qquad
\zeta = x_3.
\label{eq:indicator_target_coordinates}
\end{equation*}
By axisymmetry, we construct the tables at representative
targets $(\rho,0,\zeta)$. The unit-density indicators are
defined with respect to the coordinate directions in this
reference frame.

For each upsampling factor $\kappa$, let
$E_j^\kappa(\rho,\zeta)$ denote the unit-density indicators
$E_j$ in \eqref{eq:stresslet_collapsed_indicator}, evaluated
at $(\rho,0,\zeta)$ with the
$\kappa n_\theta\times\kappa n_\phi$ tensor-product quadrature
rule. Since these indicators vary rapidly in magnitude,
we tabulate their logarithms,
\begin{equation}
L_j^\kappa(\rho,\zeta)
\coloneqq \log_{10}E_j^\kappa(\rho,\zeta),
\qquad j=1,2,3.
\label{eq:indicator_log_table}
\end{equation}
At arbitrary targets, we evaluate bilinear interpolants of
the tabulated logarithms, denoted by the same symbols
$L_j^\kappa$.

\begin{figure}[!t]
\centering
\includegraphics[width=0.5\linewidth]{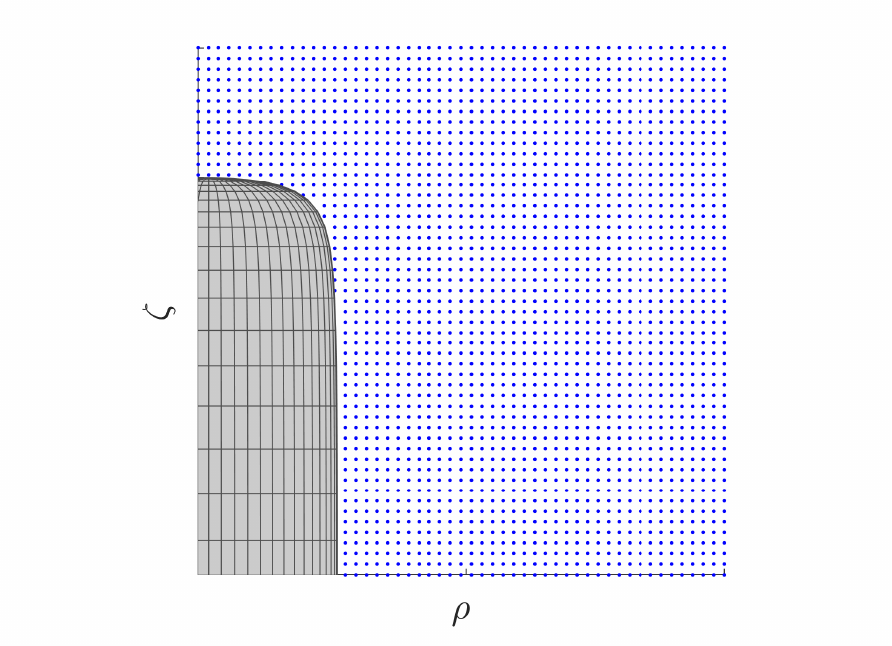}
\caption{Tabulation grid in reduced target coordinates
$(\rho,\zeta)$. Reflection symmetry of the illustrated
geometry about $\zeta=0$ allows the grid to be restricted
to $\rho\geq0$ and $\zeta\geq0$.}
\label{fig:tabulation_grid}
\end{figure}

The tables are constructed on a Cartesian grid in $(\rho,\zeta)$ that extends sufficiently far from the particle surface to capture the decay of the quadrature error with target distance. To ensure the tabulation grid is sufficiently large, we use a square grid with side length $2\times\max_{\conttheta \in [0, \pi]}(a(\conttheta) + c(\conttheta))$. Since $\rho\geq0$, only a half-plane is required. For geometries with additional reflection symmetry about $\zeta=0$, the unit-density indicators are also even in $\zeta$, and the grid can further be reduced to $\zeta\geq0$, as illustrated in Figure~\ref{fig:tabulation_grid}. For targets outside the tabulated region, we use linear extrapolation of the tabulated values. Such targets are generally sufficiently far from the surface that the corresponding quadrature errors are small. Since our primary concern is target points within the tabulated region, we use interpolation below to refer collectively to interpolation and extrapolation.

At interior tabulation points, where classification is not
required, we assign large indicator values to avoid artificially
small contributions from these points during interpolation.
At exterior points on the symmetry axis, where $\rho=0$,
the azimuthal indicator formula is undefined. We set
$E_{ij}^{\contphi}=0$ at these points and retain only the
polar contribution.

To evaluate the density modifier in
\eqref{eq:stresslet_density_modifier}, we approximate the
reference surface point by the closest node of the base grid,
that is, the original $n_\theta\times n_\phi$ surface grid.
Its parameters $(\bar{\theta},\bar{\phi})$ determine the
coordinate curves used to find the roots
$\bar{\phi}_0=\contphi_0(\bar{\theta})$ and
$\bar{\theta}_0=\conttheta_0(\bar{\phi})$.

When an explicit polar-root formula is unavailable, we
tabulate the polar root in the representative meridional
plane on the same $(\rho,\zeta)$-grid as the unit-density
indicators. We choose the root with positive imaginary part
and construct separate bilinear interpolants for its real
and imaginary parts.

We then evaluate local one-dimensional linear interpolants
of the density at $(\bar{\theta},\bar{\phi}_0)$ and
$(\bar{\theta}_0,\bar{\phi})$.
The resulting interpolated density vectors must be expressed in the same reference
frame as the tabulated unit-density indicators.
For a target with azimuth
$\alpha=\operatorname{atan2}(x_2,x_1)$,
let $B_\alpha$ denote rotation by $-\alpha$ about the symmetry
axis, mapping the target to $(\rho,0,\zeta)$. For targets on the symmetry axis, we set $\alpha=0$.
The density modifier in the reference frame is then
\begin{equation}
\sigma_j^M(\xtil)
=
\max\left\{
\left|
\left(B_\alpha\vec{\sigma}
(\bar{\theta},\bar{\phi}_0)\right)_j
\right|,
\left|
\left(B_\alpha\vec{\sigma}
(\bar{\theta}_0,\bar{\phi})\right)_j
\right|
\right\},
\qquad j=1,2,3.
\label{eq:meridional_density_modifier}
\end{equation}

Combining the interpolated unit-density indicators with the
modifier computed from the locally interpolated density gives
the following approximation to the classification indicator
in \eqref{eq:stresslet_classification_indicator}:
\begin{equation}
E^\kappa[\vec{\sigma}](\xtil)
\coloneqq
\sum_{j=1}^{3}
\sigma_j^M(\xtil)\,
10^{L_j^\kappa(\rho,\zeta)}.
\label{eq:tabulated_classification_indicator}
\end{equation}
We refer to $E^\kappa[\vec{\sigma}]$ as the
\emph{fast error indicator}.
For reflection-symmetric indicator tables, we use $|\zeta|$
in the table lookup. If the polar-root table is also restricted
to $\zeta\geq0$, the root for a target with $\zeta<0$ is
obtained as $\pi-(\bar{\theta}_0)^*$, where $\bar{\theta}_0$
is the root at $(\rho,|\zeta|)$. The density modifier is computed once per target
and reused for all $\kappa$.

Given a requested quadrature tolerance $\epsilon$ together with a set of admissible upsampling factors, we select
the smallest available $\kappa$ for which
$E^\kappa[\vec{\sigma}](\xtil)\leq\epsilon$,
with $\kappa=1$ corresponding to standard quadrature.
If no available $\kappa$ satisfies the tolerance, special quadrature is used.
Algorithms~\ref{alg:precompute} and~\ref{alg:evaluate}
summarize the main steps of the precomputation and on-the-fly
target classification, respectively. The edge cases and
symmetry conventions described above are left implicit.
An implementation of the fast classifier is available in \cite{quadind}. 

\begin{algorithm}[!t]
\caption{Precompute tabulation data for target classification}
\label{alg:precompute}
\begin{algorithmic}[1]
\Require Geometry $\boldsymbol{\gamma}$; base grid $n_\theta \times n_\phi$; upsampling factors $\{\kappa_k\}_{k=1}^{\kappamax}$; tabulation grid in $(\rho,\zeta)$
\Ensure Tables $L_j^{\kappa_k}$, $j=1,2,3$; when needed, polar-root tables $\Theta_\Re$ and $\Theta_\Im$

\If{the polar root is not available in closed form}
    \State Initialize $\Theta_\Re$ and $\Theta_\Im$ to zero
    \For{each exterior tabulation point $(\rho,\zeta)$}
        \State Set the representative target $\xtil\gets(\rho,0,\zeta)$
        \State Find the closest base-grid polar node $\bar\theta$, with $\bar\phi=0$
        \State Use Newton's method to find the polar root $\bar\theta_0\in\mathbb{C}$ satisfying $R_{\bgamma}^2(\bar\theta_0,0,\xtil)=0$
        \State Store $\Theta_\Re(\rho,\zeta)\gets\Re(\bar\theta_0)$ and $\Theta_\Im(\rho,\zeta)\gets|\Im(\bar\theta_0)|$
    \EndFor
\EndIf

\State Initialize all tables $L_j^{\kappa_k}$ with conservative large values
\For{each upsampling factor $\kappa_k$}
    \For{each exterior tabulation point $(\rho,\zeta)$}
        \State Set the representative target $\xtil=(\rho,0,\zeta)$
        \State Compute the roots and auxiliary quantities needed for the \(\kappa_k n_\theta\times \kappa_k n_\phi\) rule
        \State Evaluate directional indicators $E_{ij}^{\contphi}$ and $E_{ij}^{\conttheta}$, $i,j=1,2,3$, using \cref{eq:stresslet_unit_indicator_phi,eq:stresslet_unit_indicator_theta} and the conventions described in Section \ref{ss:tabulation}
        \For{$j=1,2,3$}
            \State Set $E_j^{\kappa_k}(\rho,\zeta)\gets\max_{i=1,2,3}E_{ij}^{\contphi}+\max_{i=1,2,3}E_{ij}^{\conttheta}$\Comment{\eqref{eq:stresslet_collapsed_indicator}}
            \State Apply a small positive floor to $E_j^{\kappa_k}(\rho,\zeta)$ before taking logarithms
            \State Store $L_j^{\kappa_k}(\rho,\zeta) \gets\log_{10}E_j^{\kappa_k}(\rho,\zeta)$\Comment{\eqref{eq:indicator_log_table}}
        \EndFor
    \EndFor
\EndFor
\end{algorithmic}
\end{algorithm}

\begin{algorithm}[!t]
\caption{Fast target classifier: Evaluate fast error indicators and classify target points}
\label{alg:evaluate}
\begin{algorithmic}[1]
\Require Geometry $\boldsymbol{\gamma}$; exterior targets $\{\xtil_m\}$; vector density $\bsigma$ on the base grid; tables from Algorithm \ref{alg:precompute}; upsampling factors $\{\kappa_k\}$; tolerance $\epsilon$
\Ensure Chosen upsampling factor for each target point, or a mark that special quadrature is required
\For{each target point $\xtil_m$}
    \State Compute cylindrical coordinates \((\rho_m,\zeta_m)\) and azimuthal angle $\alpha_m$
    \State Find the closest base-grid node $(\bar\theta,\bar\phi)$ to $\xtil_m$, choosing $\bar\phi=0$ when $\rho_m=0$
    \State Compute the azimuthal root \(\bar\phi_0=\contphi_0(\bar\theta,\xtil_m)\) analytically
    \If{the polar root is available in closed form}
        \State Compute \(\bar\theta_0=\conttheta_0(\bar\phi,\xtil_m)\) analytically
    \Else
        \State Recover \(\bar\theta_0\) by bilinear interpolation of \(\Theta_\Re\) and \(\Theta_\Im\) at \((\rho_m,\zeta_m)\)
    \EndIf
    \State Evaluate the density at \((\bar\theta,\bar\phi_0)\) and \((\bar\theta_0,\bar\phi)\) by local one-dimensional linear interpolation
    \State Apply $B_{\alpha_m}$ to interpolated density vectors and form $\sigma_j^M(\xtil_m)$, $j=1,2,3$, using \eqref{eq:meridional_density_modifier}
    \State Initially mark $\xtil_m$ as requiring special quadrature
    \For{each upsampling factor $\kappa_k$, in increasing order}
        \State Evaluate the bilinear interpolants
        $L_j^{\kappa_k}(\rho_m,\zeta_m)$,
        $j=1,2,3$
        \State Set
        $E^{\kappa_k}[\vec{\sigma}](\xtil_m)
        \gets
        \sum_{j=1}^{3}
        \sigma_j^M(\xtil_m)\,
        10^{L_j^{\kappa_k}(\rho_m,\zeta_m)}$
        \Comment{\eqref{eq:tabulated_classification_indicator}}
        \If{$E^{\kappa_k}[\vec{\sigma}](\xtil_m)\leq\epsilon$}
            \State Assign upsampling factor $\kappa_k$ to $\xtil_m$
            \State \textbf{break}
        \EndIf
    \EndFor
\EndFor
\end{algorithmic}
\end{algorithm}

\subsection{Performance of the fast error indicator}\label{ss:performance_precomputed_indicators}

In Sections~\ref{ss:error_formulas_surfaces}--\ref{ss:tabulation}, we constructed the fast error indicator for the Stokes double layer potential by combining precomputed unit-density indicators with a local density modifier. We now assess how well the resulting indicator captures the magnitude and spatial variation of the quadrature error, and measure the cost of its on-the-fly evaluation. We then compare using the local density modifier with a global density scaling for a strongly varying density.

All experiments in this section were carried out in MATLAB on a 3.4 GHz quad-core Intel i7-6700 CPU. The code used to generate the results in this section is provided in \cite{quadind}. Throughout, 
timing results are reported as the median over three runs. At each target point, the measured quadrature error is defined as the Euclidean norm of the difference between the computed and reference velocities obtained by using a large upsampling factor.

\subsubsection{Accuracy and efficiency}
As a test geometry, we consider the smooth axisymmetric capsule
\begin{equation}
\bgamma(\conttheta,\contphi) = \left(\frac{R}{C}\tanh\big(C\sin(\conttheta)\big)\cos(\contphi),\, \frac{R}{C}\tanh\big(C\sin(\conttheta)\big)\sin(\contphi),\, \frac{L}{2}\cos(\conttheta)\right),\,\conttheta\in[0,\pi],\,\contphi\in[0,2\pi),
\label{eq:capsule}
\end{equation}
where $R$ sets the radial scale, $L>0$ is the total length, and $C>0$ controls how rod-like the surface is. 
In this experiment we take $R=1$, $L=6$, and $C=3$, and discretize the surface using $n_\theta=50$ and $n_\phi=50$.
Because the polar root needed for the density modifier is not available in closed form for this geometry, we recover it by interpolation from the precomputed root table described in \Cref{ss:tabulation}.

To test the indicator with a strongly varying density, we
prescribe a smooth vector-valued density with a localized
peak on the positive $x$-side of the surface and variations
spanning approximately five orders of magnitude.\footnote{%
Specifically, we use
$\boldsymbol{\sigma}
=(1+14P\widehat{x},\,1.1+10P\widehat{z},\,0.9+12P\widehat{y})$,
where
$(\widehat{x},\widehat{y},\widehat{z})
=(\sin\vartheta\cos\varphi,\,
  \sin\vartheta\sin\varphi,\,
  \cos\vartheta)$
and
$P=10^4(P_{\mathrm{full}}-P_{\min})/(P_{\max}-P_{\min})$.
Here
$P_{\mathrm{full}}
=(1-r^2)(1-2r\widehat{x}+r^2)^{-3/2}$
with $r=0.47$, and $P_{\min}$ and $P_{\max}$ are the minimum
and maximum values of $P_{\mathrm{full}}$ on the base grid.}

\begin{figure}[t]
\centering
\begin{subfigure}[t]{0.4\textwidth}
\includegraphics[trim={2.7cm 0.0cm 2.55cm 0.3cm},clip,width=\linewidth]{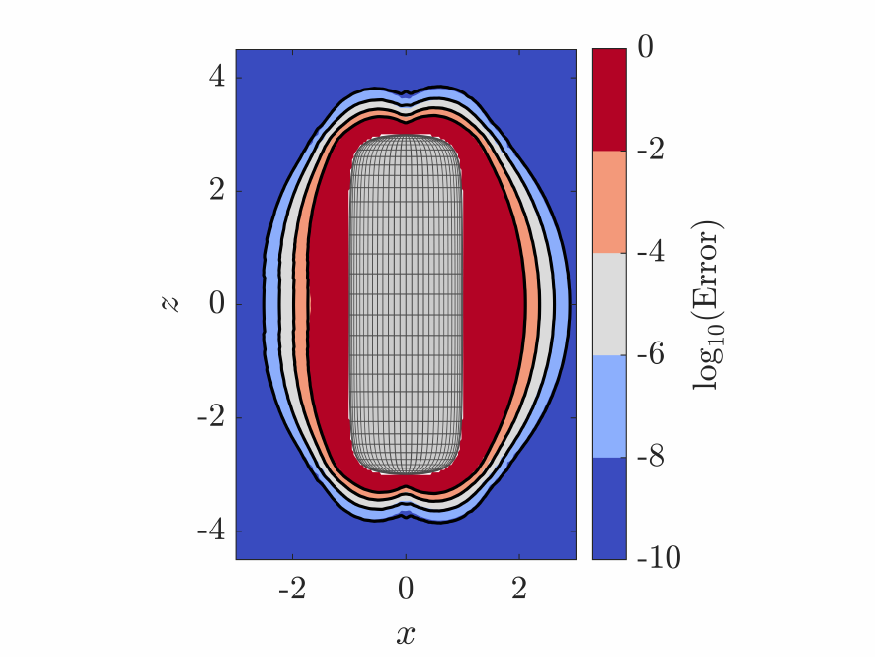}
\caption{}
\label{fig:contour_tabulated}
\end{subfigure}
\hfill
\begin{subfigure}[t]{0.49\textwidth}
\includegraphics[width=\linewidth]{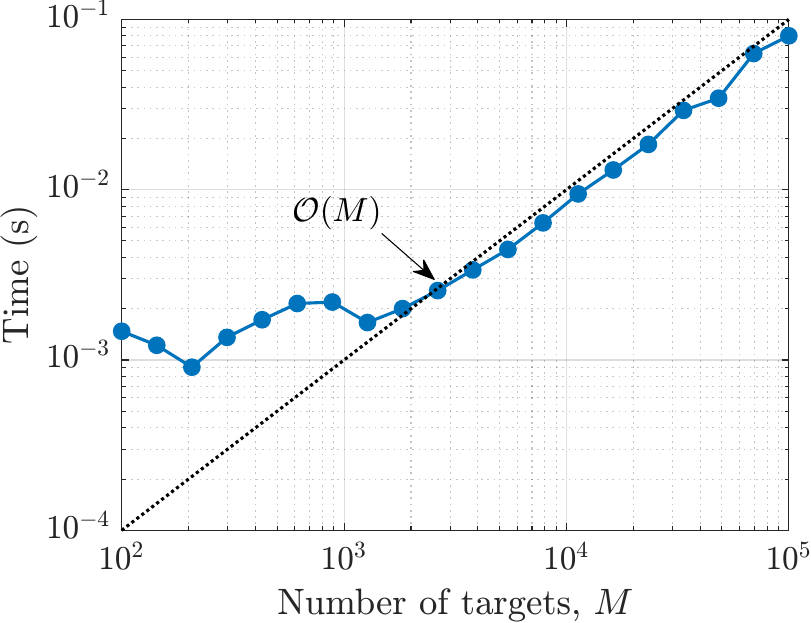}
\caption{}
\label{fig:timing}
\end{subfigure}
\caption{Panel (a) shows the measured quadrature error for the Stokes double layer potential \eqref{eq:dbl} on the capsule geometry \eqref{eq:capsule}, using
the standard quadrature rule $(\kappa=1)$. Solid black contours show the fast error indicator $E^1[\vec{\sigma}]$ at the levels $10^{-10}$, $10^{-8}$, $10^{-6}$, $10^{-4}$, and $10^{-2}$. Panel (b) shows the on-the-fly evaluation time of the fast error indicator as a function of the number of target points $M$.}
\label{fig:tabulated_vs_direct}
\end{figure}

We evaluate the indicator at varying numbers $M$ of target points in an $xz$-plane, with $y=0$, exterior to the capsule. The results are shown in Figure \ref{fig:tabulated_vs_direct}. Figure \ref{fig:contour_tabulated} compares the fast error indicator with the measured quadrature error. Despite the strongly localized density variation, the indicated error contours closely follow the measured error levels. 

Figure \ref{fig:timing} shows the on-the-fly evaluation time as a function of $M$. After a small fixed overhead for low target counts, the cost grows linearly with \(M\). For \(M=10^5\), all indicators are evaluated in less than \(0.1\) seconds, corresponding to a throughput exceeding \(10^6\) target points per second. Thus, the cost of evaluating the fast error indicator is negligible compared with the upsampled or special quadrature that it selectively activates.

\subsubsection{Effect of the density modifier}
To illustrate the importance of incorporating local information about the layer density, we compare the proposed local density modifier with a global scaling based on the $L^\infty$-norm of the density in a setting where its magnitude varies strongly over the surface. Such variation occurs, for example, in Stokes flow when two particles come very close under relative motion: in the narrow gap, the layer density becomes sharply peaked due to lubrication effects \cite{Lefebvre-Lepot_Merlet_Nguyen_2015}.

We consider the close-to-touching setup in Section \ref{ss:close_to_touching}. It consists of two prolate spheroids with equatorial radius $a=0.1$ and polar radius $c=0.5$, aligned end to end with a gap of width $2.5\times 10^{-3}$. Each particle is subject to an external force along the common polar axis, with the forces directed toward each other. This forcing drives the particles together and produces a sharply varying layer density in the narrow gap. We first solve the boundary integral equation for this density and then use it in a post-processing experiment. 

\begin{figure}[t]
\centering
\includegraphics[trim={1.5cm 0.5cm 1.5cm 0.5cm},clip,width=0.9\linewidth]{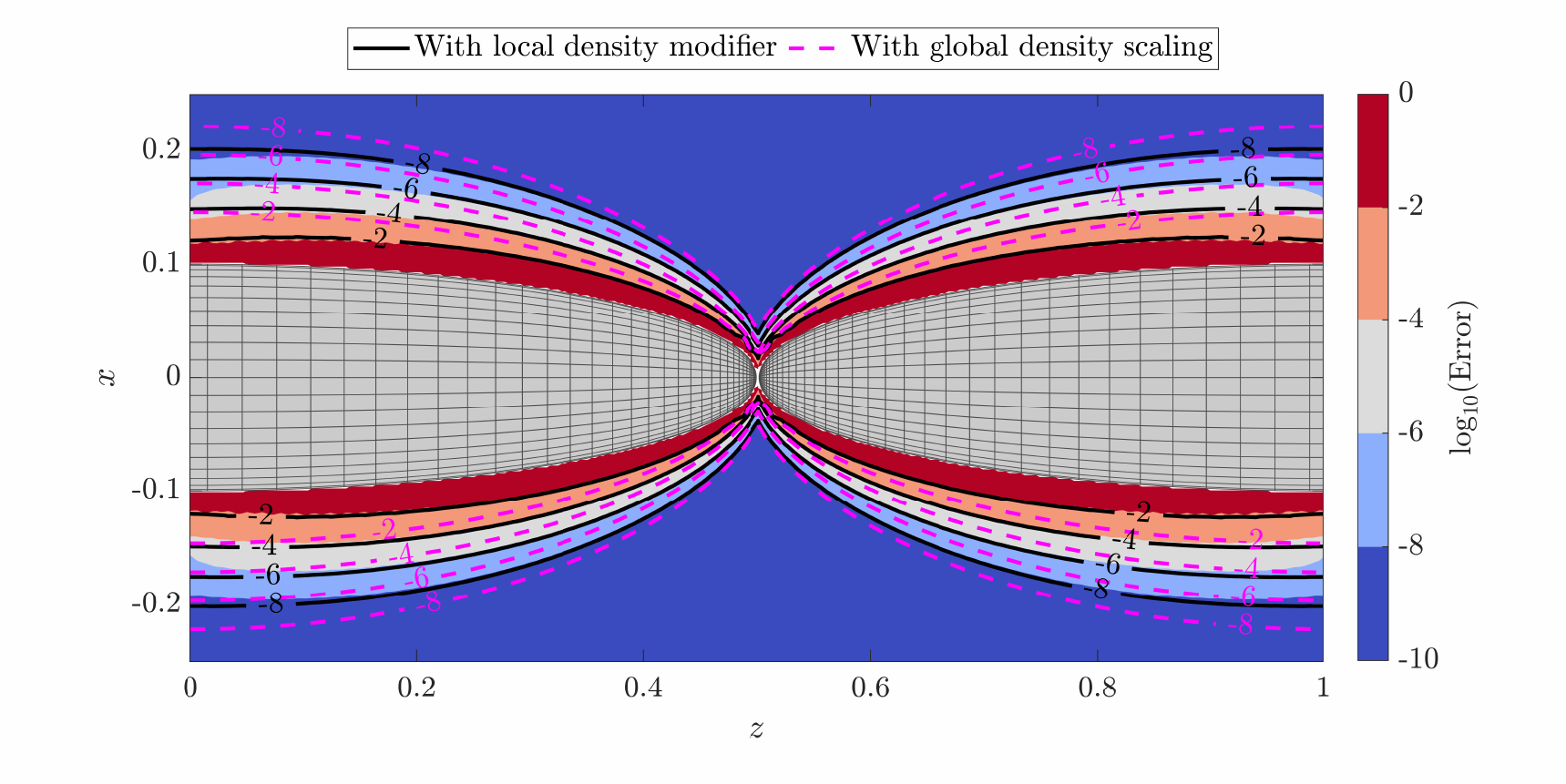}
\caption{Measured quadrature error for the Stokes double layer potential \eqref{eq:dbl} for two close-to-touching spheroids. Solid black contours show levels of the fast error indicator $E^1[\vec{\sigma}]$, which includes the local density modifier, while dashed magenta contours show the same indicator with the local modifier replaced by a global $L^\infty$-based density scaling.}
\label{fig:density_modifier}
\end{figure}

The results are shown in Figure~\ref{fig:density_modifier}. The colors show the measured quadrature error, while the solid black contours show the corresponding indicator with the density modifier included. These contours closely follow the measured error levels, demonstrating that the proposed approximation captures both the magnitude and spatial variation of the quadrature error.

The dashed magenta contours show the indicator obtained with global density scaling. They lie considerably farther from the particle surfaces, corresponding to an overestimate of approximately two orders of magnitude over much of the region shown. Although conservative, such an indicator would unnecessarily enlarge the regions assigned to upsampled or special quadrature. Including the local density modifier therefore gives a substantially sharper classification.

\section{The stabilized S3Q method}\label{s:close_evaluation}

When the target classification in \Cref{s:precomputed_error_indicators} identifies a target as requiring special quadrature, we evaluate the Stokes double layer potential \eqref{eq:dbl} using a stabilized version of S3Q \cite{krantz2026s3q}.

The present method differs from the original S3Q formulation in three main respects. 
First, the azimuthal evaluation is stabilized using TSSQ \cite{krantz2025stabilizingsingularityswapquadrature}, which addresses the severe cancellation that otherwise occurs for the stresslet kernel for very close evaluations. Second, the target-specific azimuthal SSQ weights are constructed efficiently through precomputation and one-dimensional interpolation. Third, SSQ is applied only in the azimuthal direction, while the remaining polar integral is evaluated using composite Gauss--Legendre quadrature on an adaptively refined panel grid. The recursive subdivision follows the general structure of \cite[Section~4]{krantz2026s3q}, but uses a different panel error indicator to determine whether further bisection is required.

S3Q can also be applied to a base surface discretization consisting of multiple polar panels, allowing it to be used only on panels near the target. The present implementation, however, uses the global polar discretization required by the current QBX implementation and therefore does not exploit this localization.

We begin with an overview of the stabilized method, including the density subtraction used for close targets, and then describe the azimuthal quadrature, the efficient construction of its weights, and the adaptive polar quadrature.

\subsection{Overview of the close-evaluation procedure}\label{ss:overview_close_eval}

For each target--particle interaction assigned to stabilized S3Q, we first apply density subtraction and then evaluate the resulting double layer potential as an iterated integral.

For a given exterior target point $\x$, let $\bar\x$ denote the closest node on the base surface grid of $\Gamma_q$.
Adding and subtracting the constant density $\bsigma(\bar{\x})$ gives
\begin{equation*}
\mathcal{D}_i[\bsigma;\Gamma_q](\x)
=
\mathcal{D}_i[\bsigma-\bsigma(\bar{\x});\Gamma_q](\x)
+
\mathcal{D}_i[\bsigma(\bar{\x});\Gamma_q](\x),
\qquad i=1,2,3.
\end{equation*}
Since $\x$ lies exterior to $\Gamma_q$, the second term vanishes by the stresslet identity.\footnote{For a constant density $\vec{C}$ on $\Gamma_q$, the stresslet identity reads $\bD[\vec{C};\Gamma_q](\x)=8\pi\chi_q(\x)\vec{C}$, where $\chi_q(\x)$ equals $0$, $1/2$, or $1$ when $\x$ lies exterior to, on, or interior to $\Gamma_q$, respectively \cite{pozrikidis1992boundary}.} Hence,
\begin{equation*}
\mathcal D_i[\bsigma;\Gamma_q](\x)
=
\int_{\Gamma_q}
\T_{ijk}(\x-\y)
\bigl(\sigma_j(\y)-\sigma_j(\bar{\x})\bigr)
\widehat{n}_k(\y)\,\mathrm{d}S(\y).
\label{eq:density_subtraction}
\end{equation*}
As $\x$ approaches $\Gamma_q$, the stresslet kernel develops a sharp local peak near $\bar{\x}$.
Since the subtracted density vanishes at $\bar{\x}$, this reformulation reduces the magnitude of the local peak while leaving the exact exterior potential unchanged. For notational simplicity, we continue to write $\bsigma$ for the subtracted density.

Using the surface parametrization \eqref{eq:axisymmetric_param}, we now write the Stokes double layer potential \eqref{eq:dbl} as an iterated integral, with the azimuthal integral evaluated first. Suppressing the particle index, we write the contribution from one particle surface, similarly to \eqref{eq:dbl_sum}-\eqref{eq:fij_r}, as
\begin{equation}
\mathcal{D}_i[\bsigma](\x) = \int_{0}^{\pi} U_i^{\contphi}(\conttheta;\x)\dif\conttheta,\qquad i=1,2,3,
\label{eq:dbl_close_eval}
\end{equation}
where
\begin{equation}
U_{i}^{\contphi}(\conttheta;\x) = \int_{0}^{2\pi} \frac{g_i(\conttheta,\contphi,\x)}{\big(R_{\bgamma}^2(\conttheta,\contphi,\x)\big)^{5/2}}\dif\contphi,
\label{eq:phi_integral}
\end{equation}
with
\begin{equation*}
g_i(\conttheta,\contphi,\x) = -6r_i(\vec r\cdot\widehat{\n})(\vec r\cdot\vec \bsigma)\mathcal{J}(\conttheta,\contphi),\qquad \vec r=\x-\bgamma(\conttheta,\contphi).
\label{eq:g_close_eval}
\end{equation*}
Here, $\mathcal{J}$ is the surface Jacobian defined in Section \ref{ss:axisymmetric_discretization}.

The polar interval $[0,\pi]$ is divided into panels, each equipped with an $\nGL$-point Gauss--Legendre rule. Panels are accepted or bisected using the quadrature error indicator described in Section \ref{ss:adaptiveGL}. The surface data and density values required at the polar quadrature nodes are obtained directly from the original grid by one-dimensional barycentric Lagrange interpolation.

At each polar node, the azimuthal integral \eqref{eq:phi_integral} is evaluated using the trapezoidal rule when sufficiently accurate, and otherwise using standard SSQ or its stabilized TSSQ variant, described further in Section \ref{ss:stabilized_azimuthal}. The construction of the target-specific quadrature weights is described in Section \ref{ss:efficient_evaluation_azimuthal_ssq_weights}. The resulting azimuthal integral values are then combined using composite Gauss--Legendre quadrature on the adaptively refined polar panels.

\subsection{Stabilized azimuthal evaluation}\label{ss:stabilized_azimuthal}

We begin by describing the general SSQ construction. Consider the generic line potential $I[v_p;\vec y](\xtil)$, defined on the smooth curve $\vec y$, introduced in \eqref{eq:layer_pot_curve}, repeated here for convenience,
\begin{equation}
I_p=I[v_p;\vec{y}](\xtil) = \int_\interval\frac{g(t)}{\big(R_{\vec{y}}^2(t,\xtil)\big)^p}\dif t.
\label{eq:Issq}
\end{equation}
Let $\{t_0,t_0^*\}$ denote the selected root pair of $R_{\vec{y}}^2$. The central idea in \cite{AFKLINTEBERG2021} is to transform the integral \eqref{eq:Issq} over a curve in $\mathbb{R}^3$ into an integral with known singularity structure over a simple domain in $\mathbb{C}$. This is done by introducing a simpler function $h(t,t_0)^2$ that has the same root structure as $R_{\vec y}^2(t,\xtil)$. Specifically, we then write
\begin{equation*}
I_p = \int_\interval \frac{F(t)}{h(t,t_0)^{2p}}\dif t,\qquad F(t) = g(t)\frac{h(t,t_0)^{2p}}{\big(R_{\vec y}^2(t)\big)^p}.
\end{equation*}
The ``SSQ numerator'' $F$ is smooth on $\interval$, and is approximated in a suitable basis $\{\psi_k\}_{k=1}^n$,
\begin{equation*}
F(t)\approx \sum_{k=1}^n c_k(t_0)\psi_k(t),\quad t\in \interval.
\end{equation*}
This leads to the quadrature formula
\begin{equation*}
I_p\approx\sum_{k=1}^n c_k(t_0)B_k^p(t_0),\qquad B_k^p(t_0) = \int_\interval \frac{\psi_k(t)}{h(t,t_0)^{2p}}\dif t.
\label{eq:ssq_approx}
\end{equation*}
The coefficients $c_k(t_0)$ are obtained from samples of $F$ at the quadrature nodes $\{t_j\}_{j=1}^n$, while the \emph{basis integrals} $B_k^p(t_0)$ are evaluated using formulas specific to $h$ and the chosen basis, typically through recurrence relations. The choice of $h$, basis functions $\{\psi_k\}$, and nodes $\{t_j\}$ depends on the quadrature setting in which SSQ is applied.

\subsubsection{Standard Fourier formulation}

For fixed $\conttheta$, the azimuthal integral $U_i^\contphi$ in \eqref{eq:phi_integral} is the curve integral \eqref{eq:Issq} with $\vec{y}(\contphi)=\bgamma(\conttheta,\contphi)$, $g(\contphi)=g_i(\conttheta,\contphi,\xtil)$, and $p=5/2$. We therefore apply the preceding SSQ construction to this curve. Since the procedure is identical for all the output components, we consider one component and suppress its index, together with the dependence on the fixed $\conttheta$ and $\xtil$.

Let $\contphi_0=\alpha+i\beta$, with $\alpha,\beta\in\mathbb{R}$ and $\beta>0$, denote the selected azimuthal root, computed analytically as described in Appendix \ref{a:roots}. Following \cite[Section 3]{krantz2026s3q}, we choose a Fourier basis, $h(\contphi,\contphi_0)=|e^{i\contphi}-e^{i\contphi_0}|$ and $\interval=[0,2\pi)$. For even $n_\phi$, let $c_k$ be the Fourier coefficients of $F$ at the equispaced azimuthal nodes. Standard SSQ then gives
\begin{equation*}
\label{eq:standard_fourier_expansion}
F(\varphi)
\approx
\sum_{k=-n_\phi/2}^{n_\phi/2-1}
c_k(\contphi_0)e^{ik\varphi},\quad \contphi\in[0,2\pi),
\end{equation*}
and
\begin{equation}
U^{\contphi} \approx
\sum_{k=-n_\phi/2}^{n_\phi/2-1}
c_k(\contphi_0)S_k^{5/2}(\contphi_0),
\qquad
S_k^{5/2}(\contphi_0)
=
\int_0^{2\pi}
\frac{e^{ik\contphi}}
{|e^{i\contphi}-e^{i\contphi_0}|^{5}}\dif\contphi.
\label{eq:standard_azimuthal_ssq}
\end{equation}
This is the azimuthal SSQ formulation used in the original S3Q method \cite[Section 4]{krantz2026s3q}.

\subsubsection{Stabilization using the modified Fourier basis}

In floating-point arithmetic, standard SSQ can suffer from severe cancellation when the integrand combines a strong singularity with a nearly vanishing numerator. The stresslet kernel is a particularly challenging case: for close targets, $F$ can nearly vanish around $\contphi=\alpha$, while the Fourier basis integrals $S_k^{5/2}$ become large in magnitude. The quadrature sum \eqref{eq:standard_azimuthal_ssq} may therefore involve large, nearly canceling contributions, leading to a substantial loss of accuracy \cite[Remark~11]{AFKLINTEBERG2021}, \cite[Section~3.1]{bao2024}, and \cite[Section~8.5]{krantz2026s3q}.

To stabilize the azimuthal evaluation, we use TSSQ \cite{krantz2025stabilizingsingularityswapquadrature}.
It combines a basis adapted to the local near-vanishing behavior of the numerator with stable computation of the constant coefficient. For the azimuthal integral, we use a modified Fourier basis and write:
\begin{equation*}
F(\contphi) \approx a_0(\contphi_0) + a_1(\contphi_0)\sin(\contphi-\alpha) + \sum\limits_{k=-n_\phi/2+1}^{n_\phi/2-2} b_k(\contphi_0)\sin^2\left(\dfrac{\contphi-\alpha}{2}\right)e^{ik\contphi},\quad \contphi\in[0,2\pi).
\label{eq:modified_fourier_expansion}
\end{equation*}
All nonconstant basis functions vanish at $\contphi=\alpha$. The corresponding quadrature formula is
\begin{equation}
U^{\contphi} \approx a_0(\contphi_0)S_0^{5/2}(\contphi_0) + \sum_{k=-n_\phi/2+1}^{n_\phi/2-2} b_k(\contphi_0)\widetilde{S}_k^{5/2}(\contphi_0),\qquad\widetilde{S}_k^{5/2}(\contphi_0) = \int_{0}^{2\pi} \dfrac{\sin^2\left(\frac{\contphi-\alpha}{2}\right)e^{ik\contphi}}{|e^{i\contphi}-e^{i\contphi_0}|^{5}}\dif\contphi,
\label{eq:modified_azimuthal_ssq}
\end{equation}
The integral associated with $\sin(\contphi-\alpha)$ vanishes by symmetry. We refer to \eqref{eq:modified_azimuthal_ssq} as the \emph{stabilized azimuthal approximation}.

\begin{remark}[Stable computation of the constant coefficient]\label{rem:stable_coeff_comp}
The modified basis reduces cancellation between the basis contributions, but accurate computation of the constant coefficient remains essential. For close targets, $S_0^{5/2}(\contphi_0)$ becomes large, so even a small absolute error in $a_0$ can be strongly amplified in \eqref{eq:modified_azimuthal_ssq}. Recovering $a_0$ from sampled values of $F$ can therefore lead to a loss of accuracy when $F(\alpha)$ is small. Instead, following \cite{krantz2025stabilizingsingularityswapquadrature}, we use
\begin{equation*}
a_0(\contphi_0)=F(\alpha),
\label{eq:tssq_constant_coefficient}
\end{equation*}
evaluating the known kernel and geometric factors analytically and obtaining only the density at $\contphi=\alpha$ by trigonometric interpolation from the discrete azimuthal data. The remaining coefficients $a_1$ and $b_k$ are obtained from the azimuthal samples of $F$ through the FFT-based transformation described in \cite[Section~3.2]{krantz2025stabilizingsingularityswapquadrature}.
\end{remark}

\begin{remark}[Selection of the azimuthal quadrature]
At each polar node, we first use the trapezoidal-rule indicator from Section~\ref{ss:error_formulas_one_dimension} to determine whether an azimuthal correction is required. If it indicates insufficient accuracy, we choose between standard SSQ and TSSQ using the cancellation estimate in \cite[Remark~2.1]{krantz2025stabilizingsingularityswapquadrature}. The quadrature-error indicator and cancellation estimate are evaluated separately for each output component, and their respective maxima determine a common quadrature rule for all three components.
\end{remark}

\subsection{Efficient evaluation of azimuthal SSQ weights}\label{ss:efficient_evaluation_azimuthal_ssq_weights}

The azimuthal approximations \eqref{eq:standard_azimuthal_ssq} and \eqref{eq:modified_azimuthal_ssq} combine basis coefficients with the corresponding basis integrals. In practice, we use the so-called \emph{adjoint method} to express these approximations in terms of target-specific quadrature weights, as described in \cite[Section 4.4]{krantz2026s3q}. Once these weights have been constructed for a given polar node and target point, the azimuthal integral is evaluated directly from sampled data, without explicitly forming the basis coefficients of $F$. The required density values are obtained by the polar interpolation described in Section \ref{ss:overview_close_eval}. For the modified Fourier basis, we use the stabilized adjoint formulation in \cite[Section 3.3.2]{krantz2025stabilizingsingularityswapquadrature}, together with the FFT-acceleration described in \cite[Remark 3.10]{krantz2025stabilizingsingularityswapquadrature}.

The costly part of constructing these weights is the evaluation of the standard and modified Fourier basis integrals, using the recurrence relations in \cite[Lemma 3.1]{krantz2026s3q} and \cite[Lemma 3.7]{krantz2025stabilizingsingularityswapquadrature}, respectively. For the selected root $\contphi_0=\alpha+i\beta$, with $\beta>0$, define $\chi=e^{-\beta}\in(0,1)$ and
\begin{equation*}
\mu_k^p(\chi)
\coloneqq
(1-\chi)^{2p-1}
\int_0^\pi
\frac{\cos(kt)}
{\bigl(1-2\chi\cos(t)+\chi^2\bigr)^p}
\dif t,
\qquad p=3/2,\,5/2.
\label{eq:scaled_fourier_basis_integrals}
\end{equation*}
The standard Fourier basis integrals satisfy
\begin{equation*}
S_k^{5/2}(\contphi_0)
=
\frac{2e^{ik\alpha}}{(1-\chi)^4}\mu_k^{5/2}(\chi),
\label{eq:fourier_basis_integral_scaling}
\end{equation*}
while the modified basis integrals $\widetilde{S}_k^{5/2}$ involve both $\mu_k^{5/2}$ and $\mu_k^{3/2}$. 
Thus, evaluating these two families is sufficient to construct both sets of basis integrals, with the dependence on $\alpha$ entering only through simple scaling factors.

\begin{figure}[!t]
\centering
\includegraphics[width=0.5\linewidth]{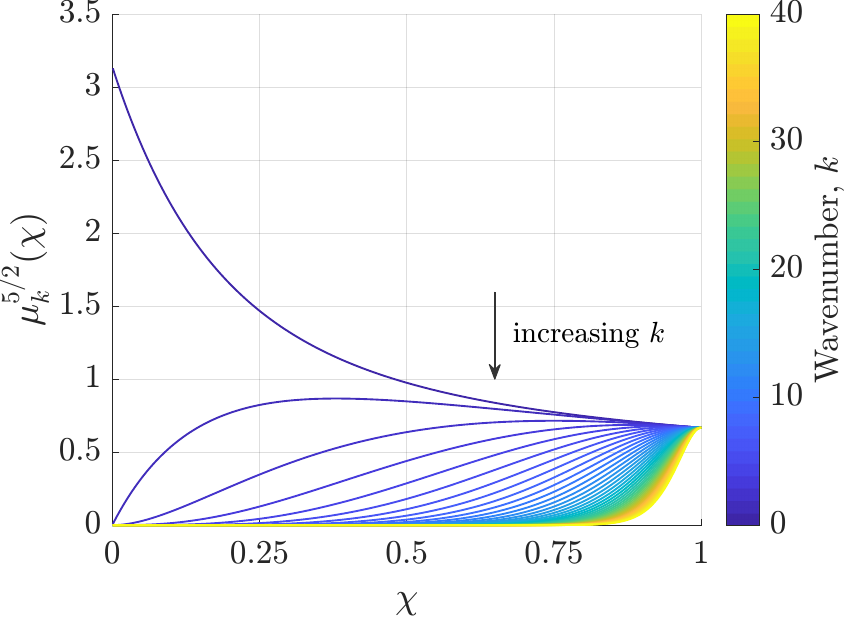}
\vspace{0.25em}
\caption{The quantities $\mu_k^{5/2}(\chi)$ as functions of $\chi$. Their smooth dependence on $\chi$ allows the azimuthal SSQ basis integrals to be evaluated efficiently to high accuracy by one-dimensional interpolation.}
\label{fig:mukp}
\end{figure}

To avoid repeatedly computing the basis integrals from the recurrence relations at every polar node, we precompute Chebyshev interpolants of $\mu_k^{3/2}$ and $\mu_k^{5/2}$ and evaluate them using the second barycentric formula \cite{trefethen2019}. 
This is effective because the quantities depend smoothly on $\chi$, as illustrated in Figure \ref{fig:mukp} for $p=5/2$. Only nonnegative modes are tabulated, since $\mu_{-k}^p=\mu_k^p$. The number of interpolation nodes is selected in a one-time offline search so that the interpolation error is below $10^{-10}$ in absolute value over the required interval in $\chi$. We use the same tables for all requested quadrature tolerances, since reducing the accuracy of the interpolation provides only negligible savings in the total evaluation time.

\subsection{Adaptive polar Gauss--Legendre quadrature} \label{ss:adaptiveGL}
After the azimuthal evaluation, the remaining polar integrand in \eqref{eq:dbl_close_eval} may still vary rapidly near the polar parameter associated with the closest surface point. We resolve this variation using composite Gauss--Legendre quadrature on an adaptively constructed panel grid. 
For a target $\xtil$, let the final partition of the base interval $[0,\pi]$ consist of $\Npan$ panels $P_m=[\conttheta_a^{(m)},\conttheta_b^{(m)}]$, $m=1,\dots,\Npan$. Applying an $\nGL$-point Gauss--Legendre rule on each panel gives
\begin{equation}
\mathcal{D}_i[\bsigma](\xtil)
\approx
\sum_{m=1}^{\Npan}
\sum_{\ell=1}^{\nGL}
w_\ell^{(m)}
U_i^{\contphi}\bigl(\conttheta_\ell^{(m)};\xtil\bigr),
\qquad i=1,2,3,
\label{eq:adaptive_polar_rule}
\end{equation}
where $\conttheta_\ell^{(m)}$ and $w_\ell^{(m)}$ are the nodes and weights mapped to $P_m$. Each azimuthal integral in \eqref{eq:adaptive_polar_rule} is evaluated using the trapezoidal rule, SSQ, or TSSQ, as described in Section \ref{ss:stabilized_azimuthal}.

\paragraph{Minimum refinement depth.}
Before applying target-dependent refinement, we uniformly bisect $[0,\pi]$ to a minimum depth
\begin{equation*}
\dmin = \max\left\{0,\left\lceil\log_2\left(\frac{n_\theta}{\nGL}\right)\right\rceil\right\},
\label{eq:dmin_heuristic}
\end{equation*}
where $n_\theta$ is the number of points in the base polar grid. 
This heuristic ensures that the initial composite grid contains at least $n_\theta$ polar nodes. In practice, we have found that this provides sufficient resolution of the geometry and density represented on the base grid, before further refinement addresses the target-dependent near singularity. 
Geometry and density values at the $\nGL$ nodes of each panel are interpolated directly from the original $n_\theta$-point grid using precomputed barycentric Lagrange interpolation matrices.

\paragraph{Target-dependent refinement.}
Starting from the minimum-depth grid, we test each panel separately and recursively bisect those that do not satisfy the refinement criterion. This criterion is based on the singularity structure of the polar integrand obtained after the azimuthal integration. This integrand contains three factors that become difficult to resolve numerically as the target approaches the surface \cite[Section 4.2]{krantz2026s3q}. Specifically, it involves a square-root-type factor, a reduced-order nearly singular factor, and a factor with a logarithmic near singularity. The nearly singular factor is the most difficult of these to resolve, even if its singularity order is one less than that of the original kernel.

Quadrature error indicators are available for the square-root and reduced-order nearly singular factors, denoted by $\quadestsq(P_m,\xtil)$ and $\quadestsing(P_m,\xtil)$, respectively. The indicator $\quadestsq(P_m,\xtil)$ is computed according to \cite[Section 5.2]{krantz2026s3q} and $\quadestsing(P_m,\xtil)$ from the formulas in Section \ref{ss:error_formulas_one_dimension}. Since the logarithmic factor was found to have resolution requirements similar to those of the square-root factor, we do not introduce a separate indicator for it. 
Since these indicators concern individual singular factors, we introduce scaling factors $C^{\mathrm{sq}}$ and $C^{\mathrm{sing}}$ to account for the remaining factors in the integrand. Like the prefactors in \eqref{eq:etz_1d}, they incorporate geometry and density information at the complex polar root. Taking the maximum over the output components separately in each scaling factor, we obtain the combined panel error indicator
\begin{equation*}
\quadest(P_m,\xtil)
=
C^{\mathrm{sq}}(\xtil)\,\quadestsq(P_m,\xtil)
+
C^{\mathrm{sing}}(\xtil)\,\quadestsing(P_m,\xtil).
\label{eq:polar_panel_indicator}
\end{equation*}
The panel $P_m$ is accepted if $\quadest(P_m,\xtil)\leq\epsilon(\conttheta_b^{(m)}-\conttheta_a^{(m)})/\pi$, where the prescribed tolerance $\epsilon$ is distributed in proportion to the panel length. 
Otherwise, the panel is bisected and the same test is applied recursively to its two children. 
The combined indicator depends on the target, geometry, density, and panel order $\nGL$. Taking the maximum over the output components gives a common refinement criterion, so all three components share the same panel grid.

\paragraph{Effect of panel order.}
The panel order $\nGL$ is fixed throughout the adaptive grid. Higher-order rules typically require fewer panels, but each panel requires more azimuthal evaluations. The computational cost therefore depends on both $\Npan$ and $\nGL$.

\begin{figure}[t]
\centering
\begin{subfigure}[t]{0.49\textwidth}
\includegraphics[width=\linewidth]{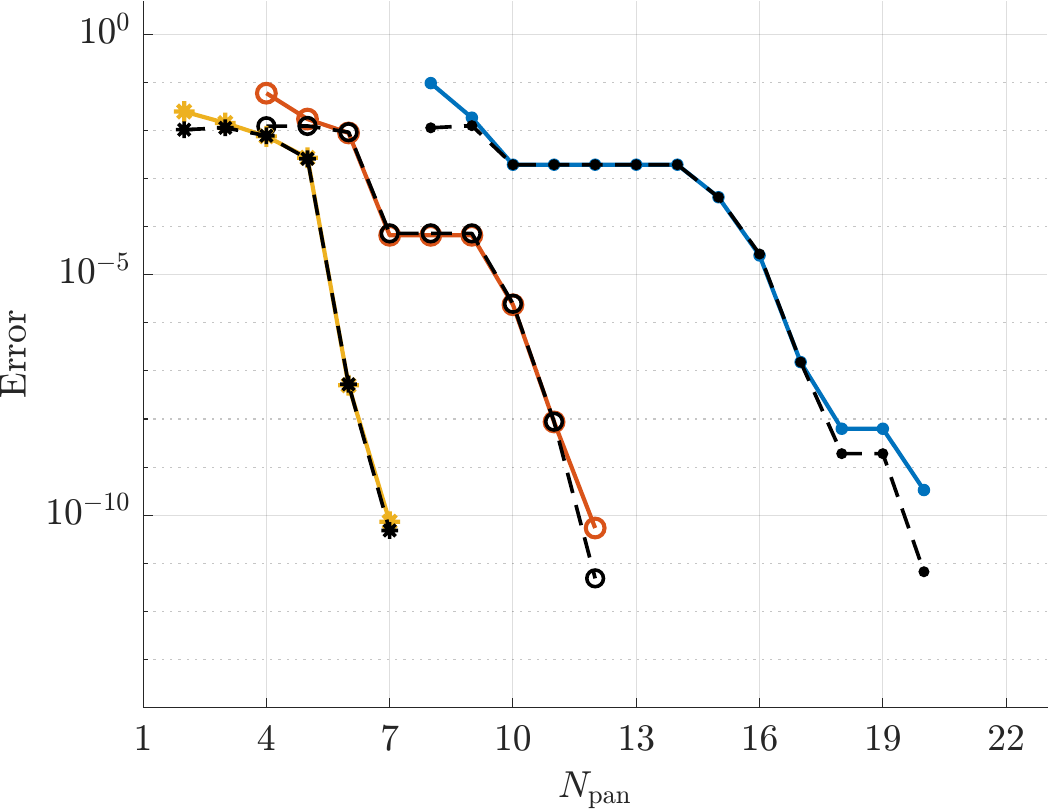}
\caption{}
\label{fig:PanelErr_1}
\end{subfigure}
\hfill
\begin{subfigure}[t]{0.49\textwidth}
\includegraphics[width=\linewidth]{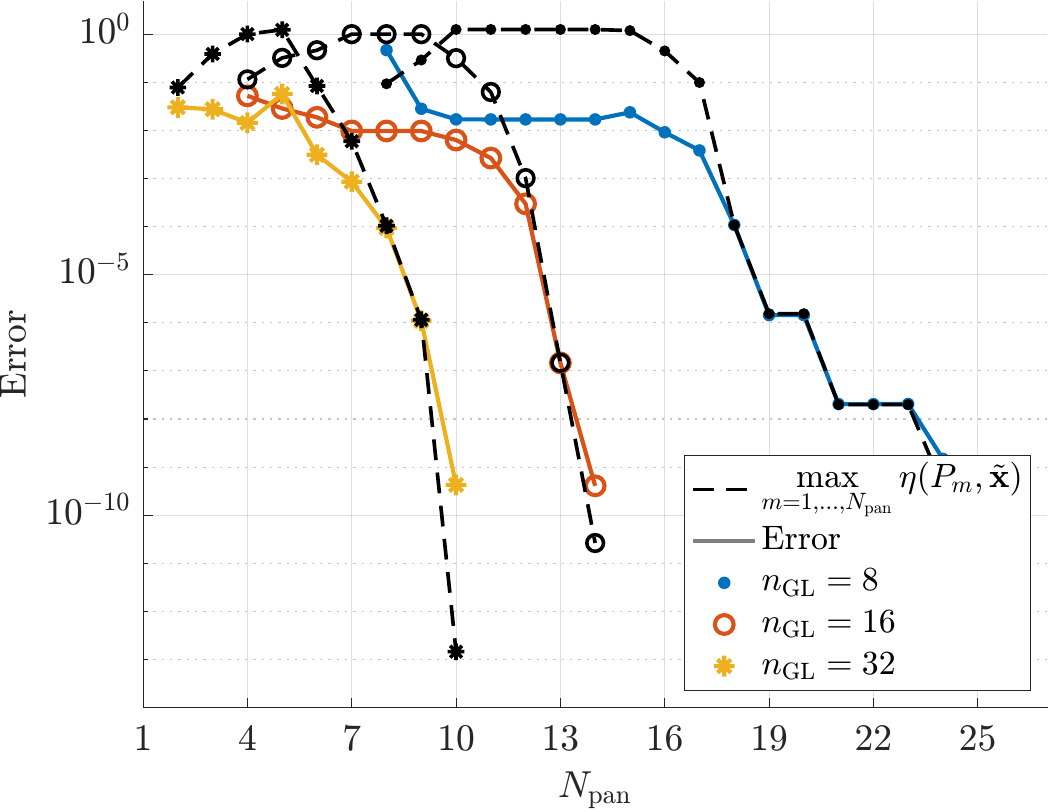}
\caption{}
\label{fig:PanelErr_2}
\end{subfigure}
\caption{Measured quadrature error of the Stokes double layer potential (solid lines) and the corresponding maximum panel error indicator (dashed lines) as functions of the number of polar panels, using \(\nGL=8\) (blue dots), \(\nGL=16\) (red open circles), and \(\nGL=32\) (yellow asterisks). Panels (a) and (b) show results for two targets on the same surface-normal line, at distances $5\times10^{-4}$ and $10^{-4}$ from the spheroid in \Cref{ss:simple}, respectively.}
\label{fig:PanelErrs}
\end{figure}

\Cref{fig:PanelErrs} compares the measured quadrature error against the maximum error indicator value over all panels for two target points along the same surface-normal line 
from the spheroid in \Cref{ss:simple}. For the target at distance $5\times 10^{-4}$ in \Cref{fig:PanelErr_1}, we see that the indicator tracks the error closely. For the closer target at distance $10^{-4}$ in \Cref{fig:PanelErr_2}, it initially overestimates the error on coarse panel grids. Moreover, with higher-order rules, a single bisection can reduce the error by several orders of magnitude, potentially leaving the achieved error far below the requested tolerance. Lower-order rules give more gradual error reductions and thus finer error control. Since all tested values of $\nGL$ can meet a specified error tolerance and have comparable computational cost, we therefore use $\nGL = 8$ in all subsequent numerical experiments.

\section{Numerical experiments}\label{s:numerical_experiments}

In \Cref{ss:performance_precomputed_indicators}, we showed that the precomputation strategy provides a fast procedure for classifying off-surface targets into regions where standard quadrature, upsampled quadrature, or special quadrature is required to meet a prescribed tolerance $\epsilon$. 
We now test the resulting tolerance-driven evaluation strategy in a more representative setting where the boundary integral equation is first solved for the layer density and the velocity field is subsequently evaluated at off-surface targets in a post-processing stage.  
Although the classifier and stabilized S3Q are formulated for general smooth axisymmetric surfaces, the experiments in this section use spheroids, for which the required on-surface QBX and off-surface QBX reference evaluations are available.

Starting with a single spheroid example in \Cref{ss:simple}, we establish the need for special quadrature beyond the simple upsampling approach and validate the use of the S3Q methodology described in \Cref{s:close_evaluation}. In \Cref{ss:cancellation_study} we then demonstrate the effect of the stabilized azimuthal evaluation (TSSQ) introduced in \Cref{ss:stabilized_azimuthal}, within S3Q. We then illustrate the effective error control and practical point categorization of the fast error indicators from \Cref{s:precomputed_error_indicators} when combined with S3Q in a two-spheroid configuration (\Cref{ss:qbx_vs_s3q}). To further demonstrate this error control, we apply this combined approach on a challenging close-to-touching particle configuration (\Cref{ss:close_to_touching}). Lastly, we consider a larger example with 100 spheroids and periodic boundary conditions, accelerated by a fast summation method (\Cref{ss:cluster}).

\subsection{Setup}
Throughout this section, we use two prolate spheroid geometries. Type~1 spheroids have equatorial radius $a=0.05$ and polar radius $c=0.1$, and are discretized using $n_\theta=40$ and $n_\phi=60$. Type~2 spheroids have $a=0.1$ and $c=0.5$, and are discretized using $n_\theta=80$ and $n_\phi=40$.

The experiments in Sections~\ref{ss:simple}--\ref{ss:close_to_touching} concern mobility problems, whereas the experiment in \Cref{ss:cluster} concerns a resistance problem. Unless stated otherwise, the boundary integral equation \eqref{eq:integral_equation} is solved using GMRES with relative residual tolerance $10^{-10}$. We use the block-diagonal single-particle preconditioner of \cite[Section~7]{bagge2021highly}, following \cite[Section~6]{AFKLINTEBERG2016}. For each spheroid type, its diagonal blocks are constructed from the explicit inverse of the corresponding single-particle system and rotated to the orientation of each particle. We use $N_{\mathrm{src}}=65$ completion sources in \eqref{eq:compflow}, placed along the polar axis inside each particle.

On-surface self-interactions are evaluated using QBX, accelerated by the axisymmetric precomputation scheme of \cite{AFKLINTEBERG2016} outlined in \Cref{ss:qbx}. The QBX parameters are initialized according to \cite[Section~6]{bagge2021highly} and then fine-tuned to ensure sufficient accuracy.

A common quadrature hierarchy is used for all off-surface target--particle interactions, both during the boundary integral solve and in post-processing. As described in \Cref{alg:evaluate}, the classifier selects standard quadrature ($\kappa=1$), the smallest adequate upsampling factor $\kappa\in\{2,\ldots,6\}$, or S3Q. The prescribed quadrature tolerance is $10^{-11}$ during the boundary integral solve and $\epsilon$ during post-processing.

In post-processing, this procedure is used to evaluate the velocity $\u(\x)$ in \eqref{eq:u_dlp} at each off-surface target $\x$. 
Reference values are computed with a tolerance of $10^{-9}$. Upsampled quadrature with $\kappa=15$ is used wherever it is predicted to meet this tolerance, and off-surface QBX, with parameters chosen for the same accuracy, is used otherwise.
As before, errors are reported as the Euclidean norm of the difference between the computed and reference velocities at each target point.

\subsection{Single spheroid and validation}\label{ss:simple}
\begin{figure}[t]
\centering
\begin{subfigure}[t]{0.43\textwidth}
\includegraphics[width=\linewidth]{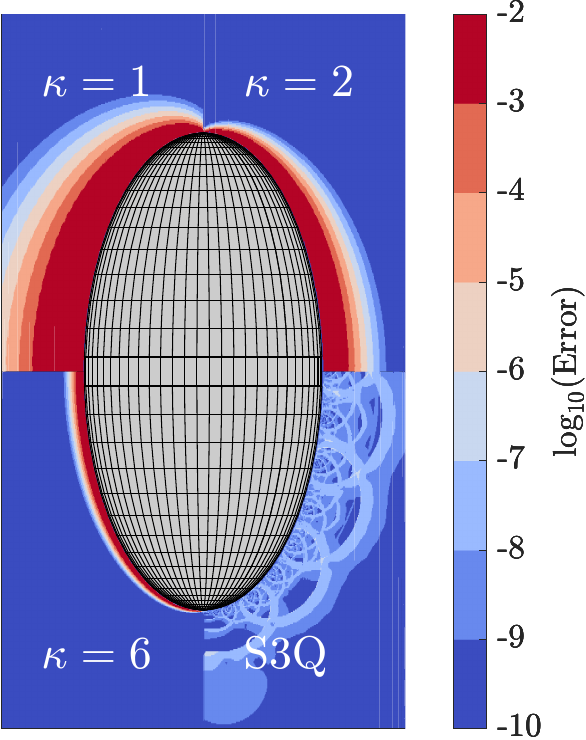}
\caption{}
\label{fig:kappa_contour}
\end{subfigure}
\hfill
\begin{subfigure}[t]{0.56\textwidth}
\includegraphics[width=\linewidth]{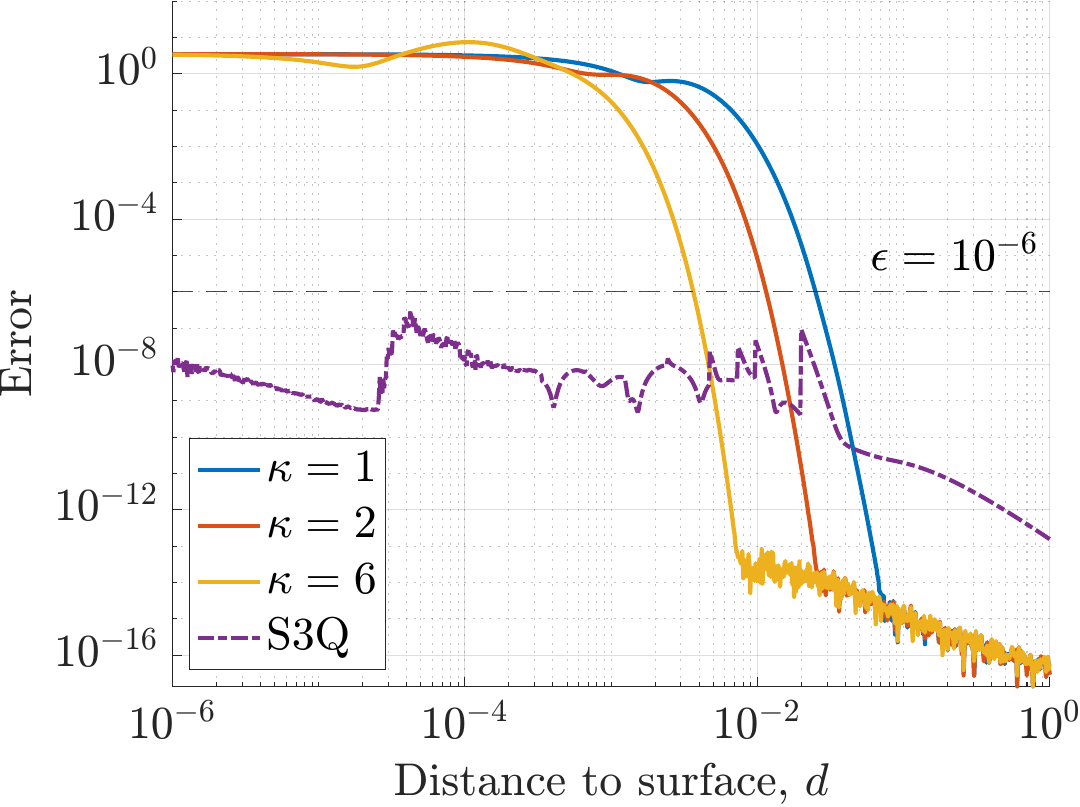}
\caption{}
\label{fig:kappa_line}
\end{subfigure}
\caption{Error comparison of standard quadrature ($\kappa = 1$), upsampled quadrature with factors $\kappa = 2$, $\kappa = 6$, and S3Q with tolerance $\epsilon=10^{-6}$. We show this for (a) a quarter plane of target points and (b) a line of target points along an outward normal from the particle surface. In panel (a), all target points were taken to be in the upper-right quarter plane for consistency, and were reflected into different regions for visualization.}
\label{fig:ExSimple}
\end{figure}

As a first test, we consider a single Type~1 spheroid subject to the (arbitrary) force $\mathbf{f}=(1,3,-2)$. We compare standard quadrature ($\kappa=1$), upsampled quadrature with $\kappa=2$ and $\kappa=6$, and S3Q with tolerance $\epsilon=10^{-6}$ on two target sets: a quarter-plane around the particle and a line of $1000$ targets extending outward along a surface normal. The resulting quadrature errors are shown in \Cref{fig:ExSimple}.

\Cref{fig:kappa_contour} shows that standard quadrature meets the prescribed tolerance for targets sufficiently far from the particle, but that its error grows rapidly as they approach the surface. Increasing $\kappa$ extends this region towards the surface, but does not eliminate the eventual error growth as the target becomes very close to the surface, even when $\kappa=6$. By contrast, S3Q maintains the error below the prescribed tolerance throughout the target set.

\Cref{fig:kappa_line} shows the same behavior over a wider range of target-to-surface distances $d$. For each fixed value of $\kappa$, the error grows and eventually plateaus as $d\rightarrow0$. For $d<10^{-3}$, the quadratures with $\kappa=1,2,6$ have lost essentially all digits of accuracy, 
whereas S3Q exhibits no such deterioration.

These results motivate the fast target classifier in \Cref{s:precomputed_error_indicators}, which invokes the more costly special quadrature (in our case S3Q) only for targets that cannot be adequately resolved using the available upsampling levels.\footnote{Notably, it is precisely in this near-surface regime that the density subtraction used in S3Q is effective; see \Cref{ss:overview_close_eval}.}

\subsection{Stabilization of S3Q near the surface}\label{ss:cancellation_study}
We now investigate the effect of the stabilized azimuthal evaluation within S3Q introduced in \Cref{ss:stabilized_azimuthal}. Here we refer to the resulting method as stabilized S3Q (using TSSQ), and to the original method, without this stabilization, as standard S3Q.
Using the Type~1 spheroid and mobility solution from the previous subsection, we place 200 random targets at a distance $d$ from the surface and evaluate the velocity field \eqref{eq:u_dlp}.
For the distance range $d\in[10^{-8},10^{-3}]$ considered here, all targets are classified as requiring special quadrature. The reported minimum, maximum, and mean errors are computed over all 200 targets at each distance.

\begin{figure}[t]
\centering
\begin{subfigure}[]{0.49\textwidth}
\includegraphics[width=\linewidth]{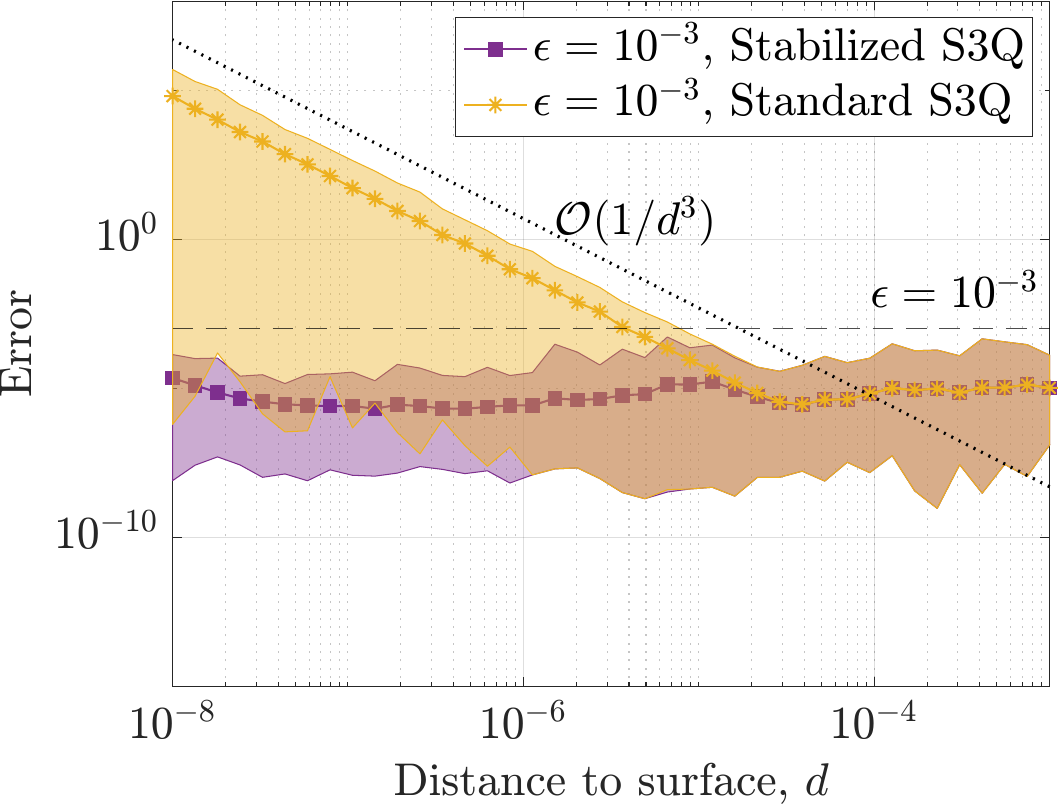}
\caption{}
\label{fig:TSSQ3}
\end{subfigure}
\begin{subfigure}[]{0.49\textwidth}
\includegraphics[width=\linewidth]{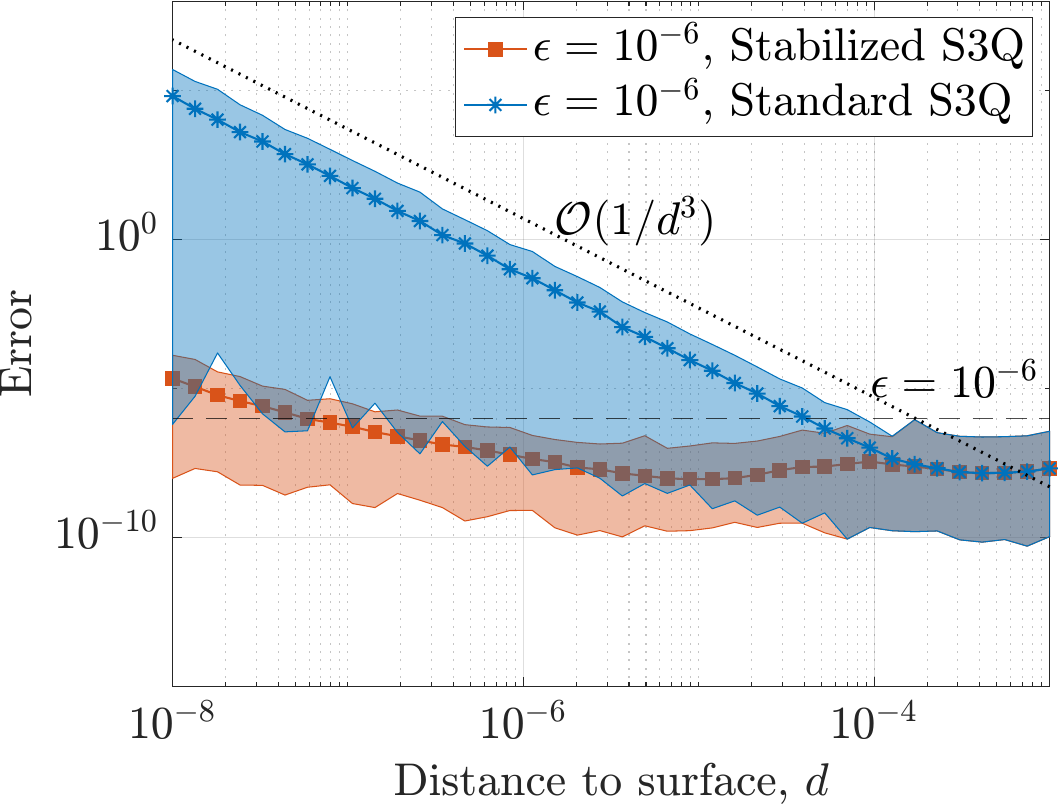}
\caption{}
\label{fig:TSSQ6}
\end{subfigure}
\caption{Panels (a) and (b) show the minimum, maximum, and mean error of the standard S3Q and stabilized S3Q methods in evaluating the velocity field \eqref{eq:u_dlp}, for tolerances $\epsilon=10^{-3}$ and $\epsilon=10^{-6}$, respectively. For each distance $d$ to the spheroid in \Cref{ss:simple}, the statistics are computed over 200 random target points.}
\label{fig:stabilization}
\end{figure}

\Cref{fig:TSSQ3,fig:TSSQ6} show the results for $\epsilon=10^{-3}$ and $\epsilon=10^{-6}$, respectively. Standard S3Q exhibits an error growth consistent with $O(1/d^3)$ independent of the requested tolerance, indicating a cancellation-dominated loss of accuracy as $d\rightarrow0$. By contrast, stabilized S3Q removes this growth and achieves errors near or below $\epsilon$ over the full range of distances. For $10^{-6}$, the maximum error is slightly above the tolerance in the closest regime, but it remains about 10 digits more accurate than the error from standard S3Q. This confirms that the loss of accuracy in standard S3Q is caused by the azimuthal SSQ step, and that replacing it by TSSQ stabilizes the close-evaluation procedure.

\begin{remark}
The $O(1/d^{3})$ growth observed for standard S3Q is consistent with the cancellation behavior reported for the Stokes double layer potential in the original S3Q paper \cite{krantz2026s3q}. In the SSQ stabilization paper \cite{krantz2025stabilizingsingularityswapquadrature}, the corresponding observed growth for a Stokes slender-body doublet kernel that has the same $|\x-\y|^{-5}$-type singularity as the stresslet was $O(1/d^{2})$. The difference is explained by the rate at which the kernel numerator vanishes near the closest point. The stresslet kernel \eqref{eq:stresslet} contains three factors of $(\x-\y)$, whereas the slender-body doublet contains two. The higher order of vanishing in the stresslet numerator leads to a stronger cancellation effect.
\end{remark}

\subsection{Evaluation of target classification and error control}\label{ss:qbx_vs_s3q}

\begin{figure}[t!]
\centering
\begin{subfigure}[t]{0.49\textwidth}
\includegraphics[width=\linewidth]{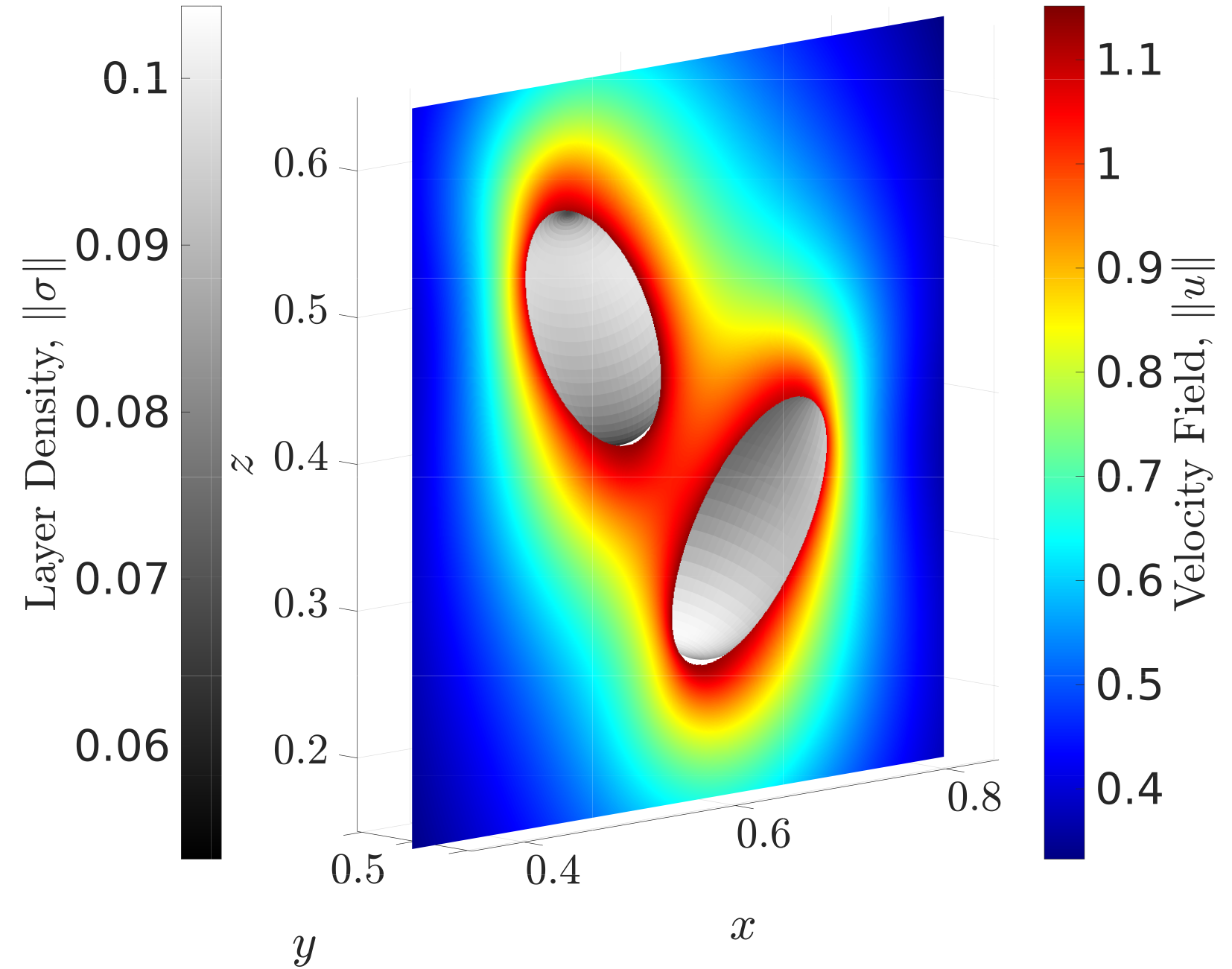}
\caption{}
\label{fig:Vel}
\end{subfigure}
\hfill
\begin{subfigure}[t]{0.49\textwidth}
\includegraphics[width=\linewidth]{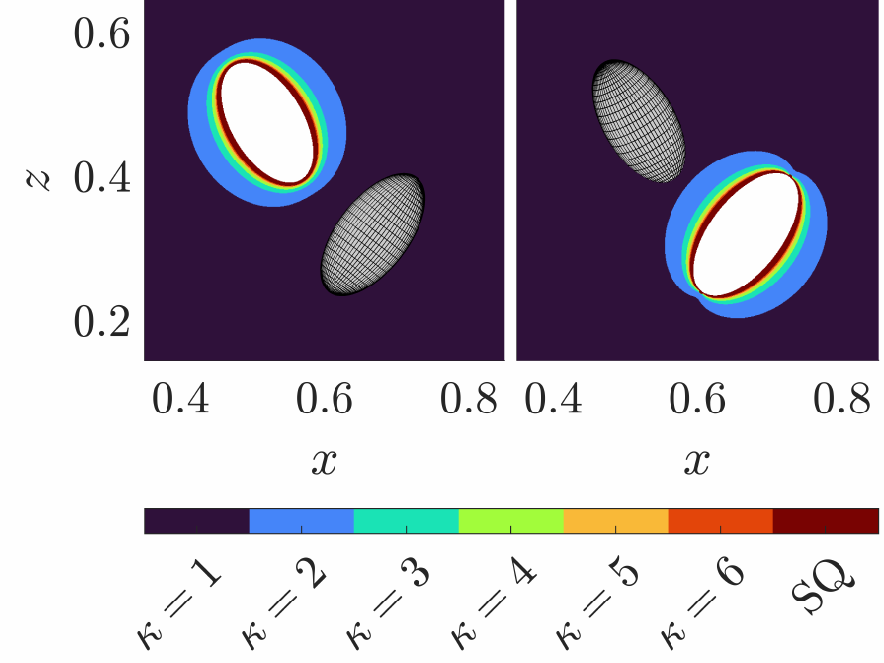}
\caption{}
\label{fig:MaskEx1}
\end{subfigure}
\caption{Panel (a) shows the velocity magnitude on the target plane together with the solved on-surface density magnitude on the spheroids. Panel (b) illustrates the resulting target point classification for $\epsilon=10^{-9}$, where $\kappa$ denotes the upsampling factor relative to the base grid.}
\label{fig:ExQBX}
\end{figure}

To investigate the efficiency of the indicator-based classification for determining the least expensive quadrature methods to obtain a prescribed accuracy, we consider two interacting particles of Type~1 with a minimum distance $7.0\times10^{-2}$ subject to the gravitational force $\mathbf{f}=(0,0,-1)$. The velocity field is evaluated on a $700\times700$ plane of target points; see \Cref{fig:Vel}. The required quadrature method for a given target point depends on the density and distance to a given particle, implying target points may require varying quadrature methods when considering the contribution from each particle, cf.~\Cref{fig:quadrature_region_illustration}. The resulting quadrature classification for each particle is shown in \Cref{fig:MaskEx1} for target points exterior to the particle. As previously expected from Figure~\ref{fig:ExSimple}, target points far from the particle surfaces can be computed simply using $\kappa=1$, and upsampled quadrature $\kappa>1$ can be leveraged to accurately evaluate at target points closer to the particle surface. For the region of points closest to the particle surface, we notice these points can only be resolved by invoking S3Q. These observations are further described in \Cref{tab:s3q_errors}, where for three tolerances $\epsilon = 10^{-3}, 10^{-6}, 10^{-9}$, the classification, mean, and max errors are summarized. The number of target points computed using each quadrature method is noted, and the total number of points for each tolerance is equal to twice the number of targets exterior to both particles.

Most targets are handled by inexpensive standard or upsampled quadrature, while the more costly special quadrature method is only invoked for a small set of targets. 
For upsampled quadrature, the cost of evaluating the quadrature sum at each target increases proportionally to $\kappa^2$, simply due to the number of points in the refined grid. The larger grid likewise  increases the cost of interpolating the geometry and layer density from the base grid, but this cost is amortized over all targets evaluated using the same $\kappa$.
For all tolerances, the throughput of S3Q, i.e.~the number of target points evaluated per second, is typically one order of magnitude lower than that of upsampled quadrature with the largest factor used, $\kappa=6$. 
As discussed in \Cref{s:close_evaluation}, the discretization used here is global and therefore does not permit the S3Q correction to be restricted to panels near the target. A multi-panel based discretization, together with further implementation optimization, could therefore improve the throughput of S3Q.

\Cref{tab:s3q_errors} shows the desired mean accuracy was met in all cases, with maximum errors typically on the same order as the tolerance, where only for $\epsilon=10^{-9}$ we notice a marginally larger factor. 
To investigate to what extent the fast target classification is optimal, we now consider how the  results change if we shift the classification in  \Cref{tab:s3q_errors} by one level. This means that for each target point, we use the next less expensive quadrature method than was initially assigned. The results are summarized in \Cref{tab:cate_compare}. For example, consider the points that were classified as needing an upsampling factor of $\kappa=4$ for the tolerance $10^{-6}$. Then only one target point out of $4447$ points has an error that is larger than the tolerance, shown as $0.02\%$ in the table. Using instead $\kappa=3$ for evaluation for the same set of points, $4246$ points, or $95.5\%$ of the points, now exceed the tolerance. 
The same pattern is observed across the other classification regions. Among all $59\,696$ points assigned to upsampling or S3Q at this tolerance, the fraction exceeding the tolerance increases from $0.02\%$ with the original assignments to more than $95\%$ after shifting them down by one level, with many of the resulting errors exceeding the tolerance by a large margin.
This example shows that, in nearly all cases, the classifier selects the least costly quadrature treatment that meets the prescribed tolerance.

\begin{table}[t]
\centering
\caption{Target classification, mean, and max errors of the quadratures used for the two-particle setup of \Cref{ss:qbx_vs_s3q}, shown in \Cref{fig:ExQBX}. Results are reported for the three tolerances $\epsilon=10^{-3},10^{-6},10^{-9}$.}
\setlength{\tabcolsep}{3pt}
\resizebox{1\textwidth}{!}{\begin{minipage}{\textwidth}
\begin{tabular}{lrrrrrrrrr}
\toprule
&
\multicolumn{3}{c}{Number of targets}
&
\multicolumn{3}{c}{Mean error}
&
\multicolumn{3}{c}{Max error}
\\
\cmidrule(lr){2-4}
\cmidrule(lr){5-7}
\cmidrule(lr){8-10}
{\bf Method}
& \multicolumn{1}{c}{$\epsilon=10^{-3}$}
& \multicolumn{1}{c}{$\epsilon=10^{-6}$}
& \multicolumn{1}{c}{$\epsilon=10^{-9}$}
& \multicolumn{1}{c}{$\epsilon=10^{-3}$}
& \multicolumn{1}{c}{$\epsilon=10^{-6}$}
& \multicolumn{1}{c}{$\epsilon=10^{-9}$}
& \multicolumn{1}{c}{$\epsilon=10^{-3}$}
& \multicolumn{1}{c}{$\epsilon=10^{-6}$}
& \multicolumn{1}{c}{$\epsilon=10^{-9}$}
\\
\midrule
Std.~quad
& 831771 & 802464 & 765355
& \(3.65\times10^{-6}\) & \(4.28\times10^{-9}\) & \(5.29\times10^{-12}\)
& \(1.11\times10^{-3}\) & \(9.63\times10^{-7}\) & \(8.75\times10^{-10}\)
\\
\(\kappa=2\)
& 15985 & 33297 & 57364
& \(9.22\times10^{-5}\) & \(4.61\times10^{-8}\) & \(2.81\times10^{-11}\)
& \(1.21\times10^{-3}\) & \(1.09\times10^{-6}\) & \(9.81\times10^{-10}\)
\\
\(\kappa=3\)
& 4913 & 9400 & 14590
& \(1.97\times10^{-4}\) & \(1.40\times10^{-7}\) & \(6.89\times10^{-11}\)
& \(1.17\times10^{-3}\) & \(1.09\times10^{-6}\) & \(9.42\times10^{-10}\)
\\
\(\kappa=4\)
& 2346 & 4447 & 6687
& \(2.90\times10^{-4}\) & \(1.63\times10^{-7}\) & \(1.11\times10^{-10}\)
& \(1.05\times10^{-3}\) & \(1.01\times10^{-6}\) & \(1.11\times10^{-9}\)
\\
\(\kappa=5\)
& 1328 & 2569 & 3814
& \(3.70\times10^{-4}\) & \(2.17\times10^{-7}\) & \(1.46\times10^{-10}\)
& \(1.24\times10^{-3}\) & \(1.15\times10^{-6}\) & \(1.05\times10^{-9}\)
\\
\(\kappa=6\)
& 856 & 1670 & 2499
& \(4.28\times10^{-4}\) & \(2.62\times10^{-7}\) & \(1.91\times10^{-10}\)
& \(1.39\times10^{-3}\) & \(9.91\times10^{-7}\) & \(1.09\times10^{-9}\)
\\
S3Q
& 4961 & 8313 & 11851
& \(4.82\times10^{-5}\) & \(1.66\times10^{-8}\) & \(1.42\times10^{-10}\)
& \(4.29\times10^{-3}\) & \(1.80\times10^{-6}\) & \(3.40\times10^{-8}\)
\\
\bottomrule
\end{tabular}
\end{minipage}}
\label{tab:s3q_errors}
\end{table}

\begin{table}[t]
\centering
\caption{Ratio of target points above the prescribed tolerances, $\epsilon=10^{-3},10^{-6},10^{-9}$, when applying the fast-target classification shown in \Cref{tab:s3q_errors} versus using the next less expensive quadrature method.}
\setlength{\tabcolsep}{3pt}
\resizebox{1\textwidth}{!}{\begin{minipage}{\textwidth}
\begin{tabular}{rrrrrrrrrrr}
\toprule
\multicolumn{3}{c}{Number of targets}
&
\multicolumn{4}{c}{Points above tolerance}
&
\multicolumn{4}{c}{Points above tolerance}
\\
\cmidrule(lr){1-3}
\cmidrule(lr){4-7}
\cmidrule(lr){8-11}
\multicolumn{1}{c}{$\epsilon=10^{-3}$}
& \multicolumn{1}{c}{$\epsilon=10^{-6}$}
& \multicolumn{1}{c}{$\epsilon=10^{-9}$}
& {\bf Method}
& \multicolumn{1}{c}{$\epsilon=10^{-3}$}
& \multicolumn{1}{c}{$\epsilon=10^{-6}$}
& \multicolumn{1}{c}{$\epsilon=10^{-9}$}
& {\bf Method}
& \multicolumn{1}{c}{$\epsilon=10^{-3}$}
& \multicolumn{1}{c}{$\epsilon=10^{-6}$}
& \multicolumn{1}{c}{$\epsilon=10^{-9}$}
\\
\midrule
831771 & 802464 & 765355
& Std.~quad
& \(<0.01 \%\) & \(0 \%\) & \(0 \%\)
& \multicolumn{1}{c}{---} 
& \multicolumn{1}{c}{---} & \multicolumn{1}{c}{---} & \multicolumn{1}{c}{---}
\\
15985 & 33297 & 57364
& \(\kappa=2\)
& \(0.01 \%\) & \(<0.01 \%\) & \(0 \%\)
& Std.~quad 
& \(93.8 \%\) & \(95.1 \%\) & \(95.3 \%\)
\\
4913 & 9400 & 14590
& \(\kappa=3\)
& \(0.1 \%\) & \(0.01 \%\) & \(0 \%\)
& \(\kappa=2\)
& \(93.3 \%\) & \(95.9 \%\) & \(96.2 \%\)
\\
2346 & 4447 & 6687
& \(\kappa=4\)
& \(0.04 \%\) & \(0.02 \%\) & \(0.04 \%\)
& \(\kappa=3\)
& \(92.3 \%\) & \(95.5 \%\) & \(96.1 \%\)
\\
1328 & 2569 & 3814
& \(\kappa=5\)
& \(0.6 \%\) & \(0.15 \%\) & \(0.02 \%\)
& \(\kappa=4\)
& \(89.5 \%\) & \(95.4 \%\) & \(96.2 \%\)
\\
856 & 1670 & 2499
& \(\kappa=6\)
& \(0.8 \%\) & \(0 \%\) & \(0.1 \%\)
& \(\kappa=5\)
& \(84.6 \%\) & \(94.3 \%\) & \(95.3 \%\)
\\
4961 & 8313 & 11851
& S3Q
& \(0.8 \%\) & \(0.03 \%\) & \(0.5 \%\)
& \(\kappa=6\)
& \(91.2 \%\) & \(98.0 \%\) & \(98.8 \%\)
\\
\bottomrule
\end{tabular}
\end{minipage}}
\label{tab:cate_compare}
\end{table}

\subsection{Close-to-touching spheroids}\label{ss:close_to_touching}
Next, we consider a more challenging configuration of two close-to-touching Type~2 spheroids. The particles are aligned with their polar axes and placed end to end, with a gap of width $2.5\times10^{-3}$ between their tips. They are driven toward each other by the forces $\mathbf{f}_1=(0,0,-1)$ and $\mathbf{f}_2=(0,0,1)$. 
Through lubrication effects, the resulting relative motion produces a sharply varying layer density near the facing tips, making this configuration a demanding test of the complete close-evaluation workflow \cite{Lefebvre-Lepot_Merlet_Nguyen_2015}.

\begin{figure}[t]
\centering
\begin{subfigure}[t]{0.49\textwidth}
\includegraphics[width=\linewidth]{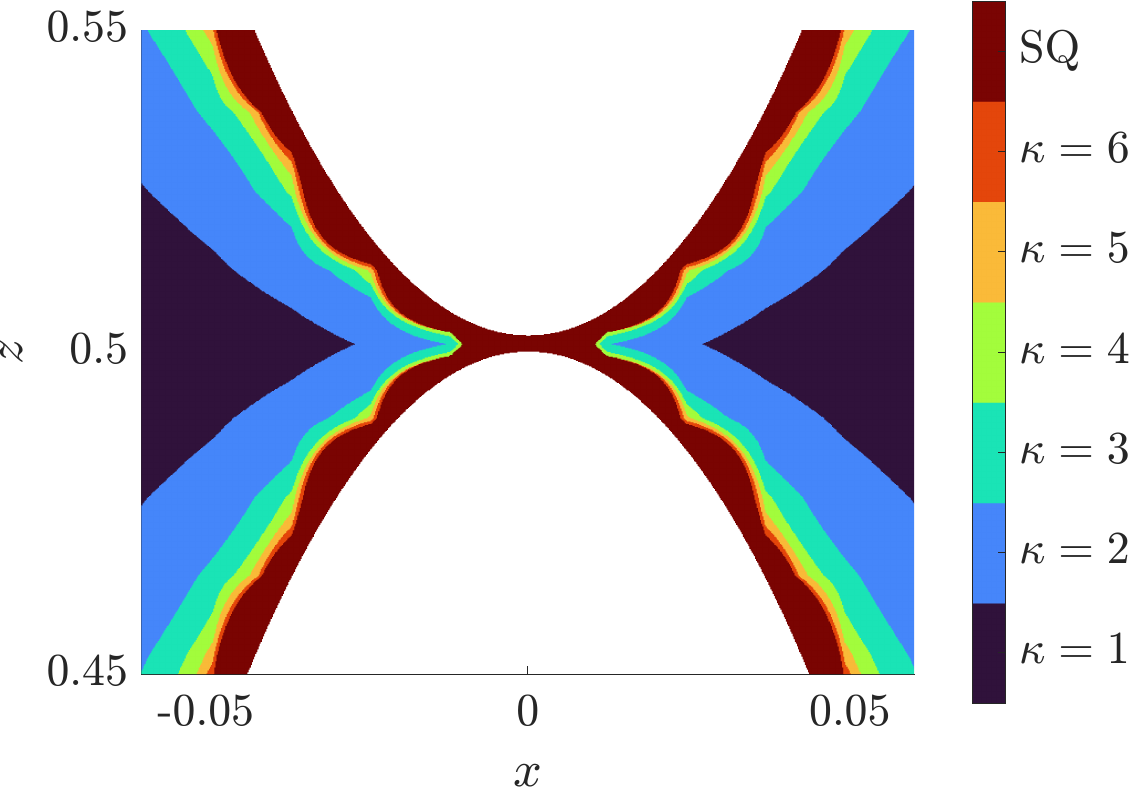}
\caption{}
\label{fig:CloseCate}
\end{subfigure}
\hfill
\begin{subfigure}[t]{0.49\textwidth}
\includegraphics[width=\linewidth]{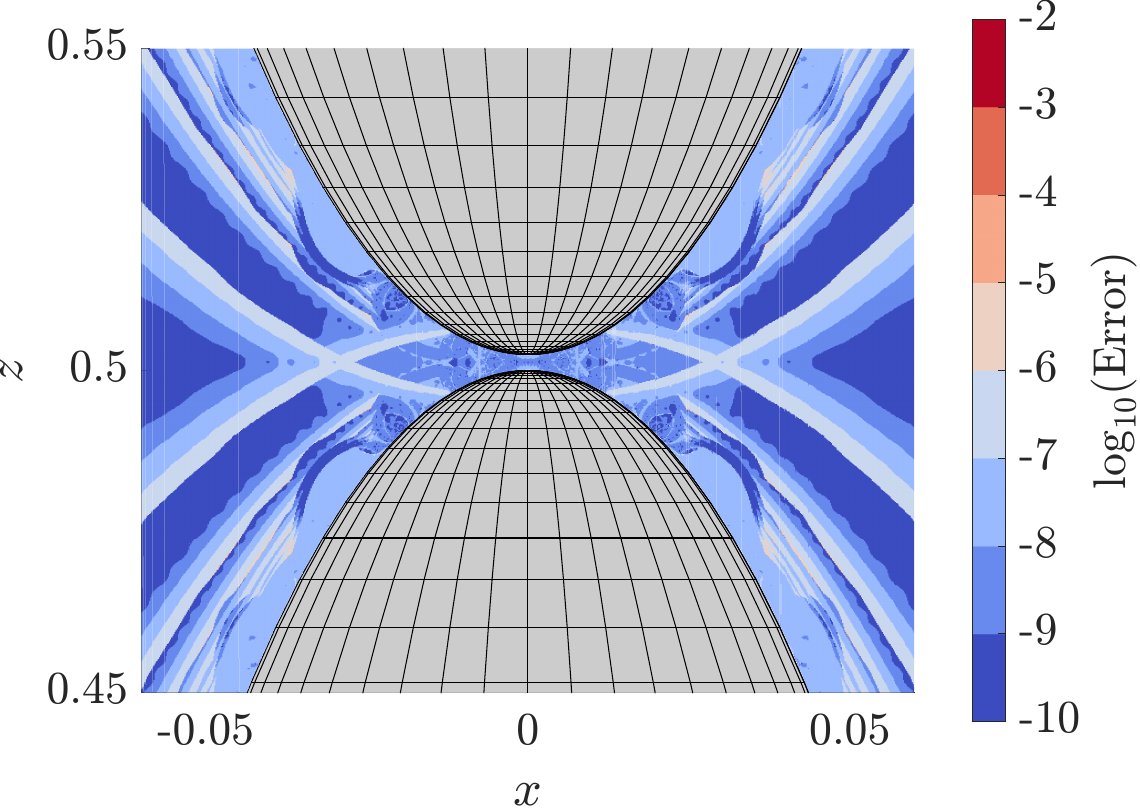}
\caption{}
\label{fig:CloseErr}
\end{subfigure}
\caption{Panel (a) shows the target classification where only the most expensive quadrature used for each point is displayed, with $\kappa$ denoting the upsampling factor relative to the base grid of the Type~2 particles. Panel (b) shows the error in evaluating the Stokes double layer potential \eqref{eq:u_dlp} using the assigned quadrature method in panel (a).}
\label{fig:ExClose}
\end{figure}

We place a $700\times 700$ target grid in the plane $y=0$ across the gap region, discard targets inside the particles, classify the remaining targets, and evaluate the velocity using the assigned quadrature method. \Cref{fig:ExClose} shows the result for $\epsilon=10^{-6}$. Classification for each target point is shown in \Cref{fig:CloseCate}, performed per particle for each target point, as showcased for the previous example in 
\Cref{fig:MaskEx1}. For visualization purposes, \Cref{fig:CloseCate} shows only the more expensive of the two quadrature methods assigned at each target. The corresponding error values are shown in \Cref{fig:CloseErr}, where we note the visible error contours produced by the transition of quadrature methods indicated in \Cref{fig:CloseCate}. The on-the-fly indicators remain effective in this close-to-touching case: S3Q is activated for the closest targets, including targets as close as $1.9 \times 10^{-5}$ to a particle surface. Measured error quantities remain close to the prescribed tolerance, with an average and maximum error of $8.46 \times 10^{-8}$ and $9.4 \times 10^{-6}$, respectively.

\subsection{Triply periodic flow through a spheroid cluster}\label{ss:cluster}
As a larger-scale demonstration, we consider 100 Type~1 spheroids placed at random positions and orientations in a rectangular cell with periodic boundary conditions in all three coordinate directions.
The configuration has a minimum separation of $2.6\times10^{-3}$, including across periodic boundaries.
Each particle has 2400 discretization points, giving a total of 720\,000 scalar unknowns. We solve the resistance problem with all particles held fixed, impose the no-slip boundary condition on their surfaces, and prescribe the background flow $\u_{\mathrm{bg}}=(1,0,0)$.\footnote{In the Fourier representations of the periodic Stokes potentials, the zero modes contain undetermined constants because the fundamental solutions are defined only up to additive constants. We adopt the usual convention of fixing these constants by requiring the periodic potentials to generate zero mean flow. For the Stokeslet and rotlet, this amounts to setting the zero mode to zero, whereas the stresslet requires an additional zero-mode term; see \cite[Sections~3.2--3.3]{AfKlinteberg2014} and \cite[Appendix~A.1]{BAGGE2023}.}

The standard quadrature sums are accelerated by the periodic DMK method \cite{krantz2026dmk}, while QBX and the proposed off-surface quadrature strategy treat the singular and nearly singular interactions, respectively.
We solve for the layer density using a relative GMRES tolerance of $10^{-8}$ and S3Q and DMK tolerances of $10^{-10}$, and allow upsampled quadrature with $\kappa=2,\ldots,6$, in addition to standard quadrature $\kappa=1$. 
The particle configuration and the resulting flow streamlines are shown in \Cref{fig:cluster_streamlines}. 

We use two error measures to assess the numerical solution. First, the normalized pointwise no-slip error is defined by
\begin{equation*}
e_{\mathrm{ns}}(\boldsymbol{x})
=
\frac{
\left\|
\boldsymbol{u}_{\mathrm{ref}}(\boldsymbol{x})
-
\boldsymbol{u}_{\Gamma}(\boldsymbol{x})
\right\|_2
}{
\left\|
\boldsymbol{u}_{\mathrm{bg}}
\right\|_2
},
\qquad
\boldsymbol{x}\in\Gamma,
\label{eq:cluster_noslip_error}
\end{equation*}
where $\boldsymbol{u}_{\mathrm{ref}}$ is the total velocity computed from the solved density with the off-surface quadrature tolerance $\epsilon=10^{-10}$, and $\boldsymbol{u}_{\Gamma}$ is the prescribed rigid-body velocity on the particle surfaces (we impose no-slip boundary conditions, $\boldsymbol{u}_{\Gamma}=\boldsymbol{0}$).

Second, to isolate the error introduced by the off-surface quadrature strategy, we define
\begin{equation*}
e_{\mathrm{quad}}(\boldsymbol{x};\epsilon)
=
\frac{
\left\|
\boldsymbol{u}_{\epsilon}(\boldsymbol{x})
-
\boldsymbol{u}_{\mathrm{ref}}(\boldsymbol{x})
\right\|_2
}{
\left\|
\boldsymbol{u}_{\mathrm{bg}}
\right\|_2
},
\qquad
\boldsymbol{x}\in\Gamma,
\label{eq:cluster_quadrature_error}
\end{equation*}
where $\boldsymbol{u}_{\epsilon}$ is obtained by evaluating the same density with tolerance $\epsilon$. 

The two errors are evaluated on a 2400-point check grid on each of five selected particle surfaces. Each check grid has the same resolution as the base discretization but is shifted in the azimuthal parameter, so that its nodes do not coincide with the collocation points used in the boundary integral solve. We write the discrete area-normalized $L^2$ error and the maximum pointwise error over the five surfaces as $\|e\|_{2}^{\mathrm{chk}}$ and $\|e\|_{\infty}^{\mathrm{chk}}$, respectively.

\begin{figure}[t!]
\centering
\begin{minipage}[t]{0.48\textwidth}
\centering
\includegraphics[width=\linewidth]{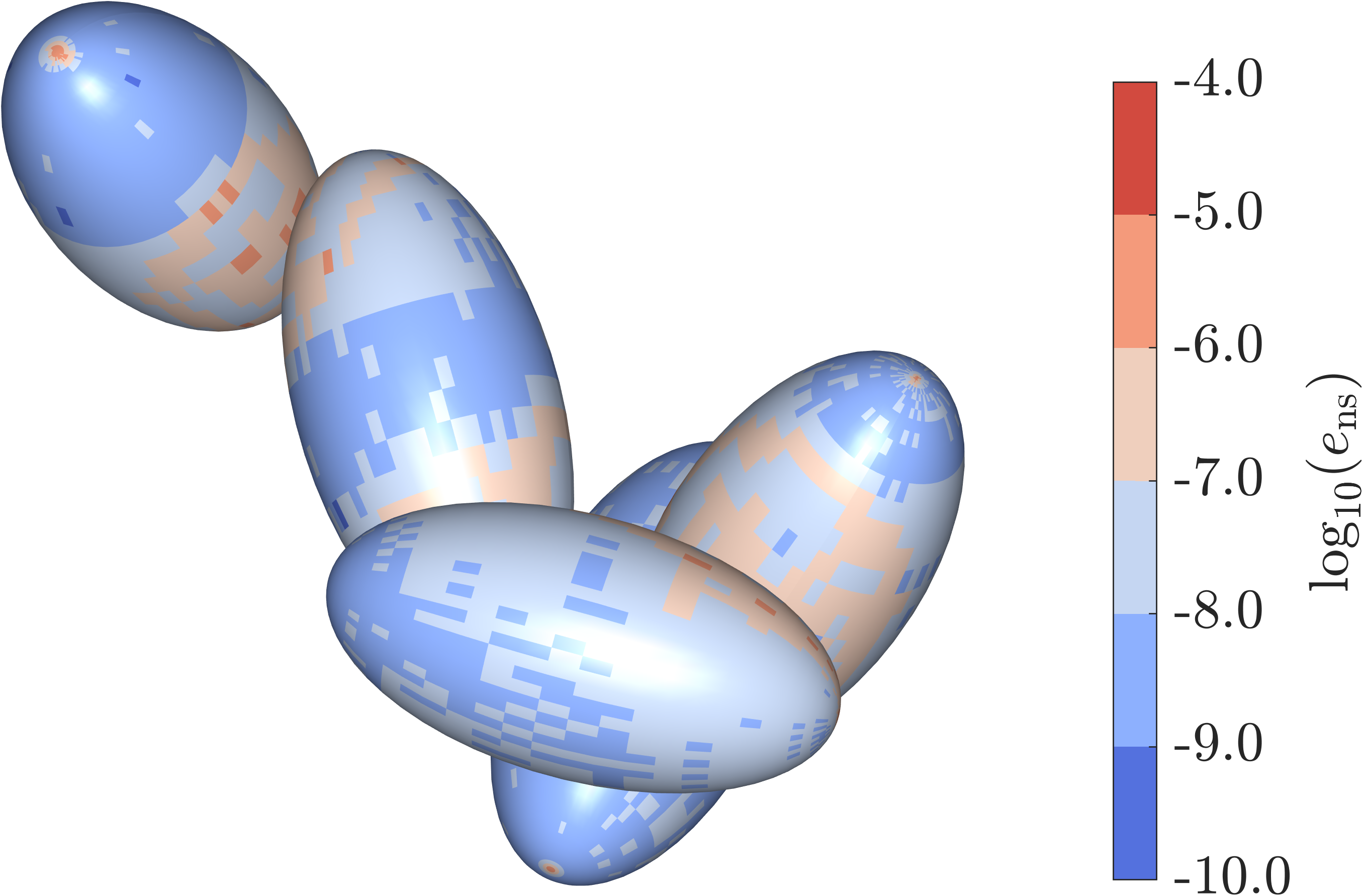}
\vspace{0.1em}
\caption{Normalized pointwise no-slip error $e_{\mathrm{ns}}$ on five Type~1 spheroids, with equatorial radius $a=0.05$ and polar radius $c=0.1$, selected from the triply periodic simulation of $100$ particles.}
\label{fig:cluster_residual}
\end{minipage}%
\hfill
\begin{minipage}[t]{0.48\textwidth}
\centering
\includegraphics[width=\linewidth]{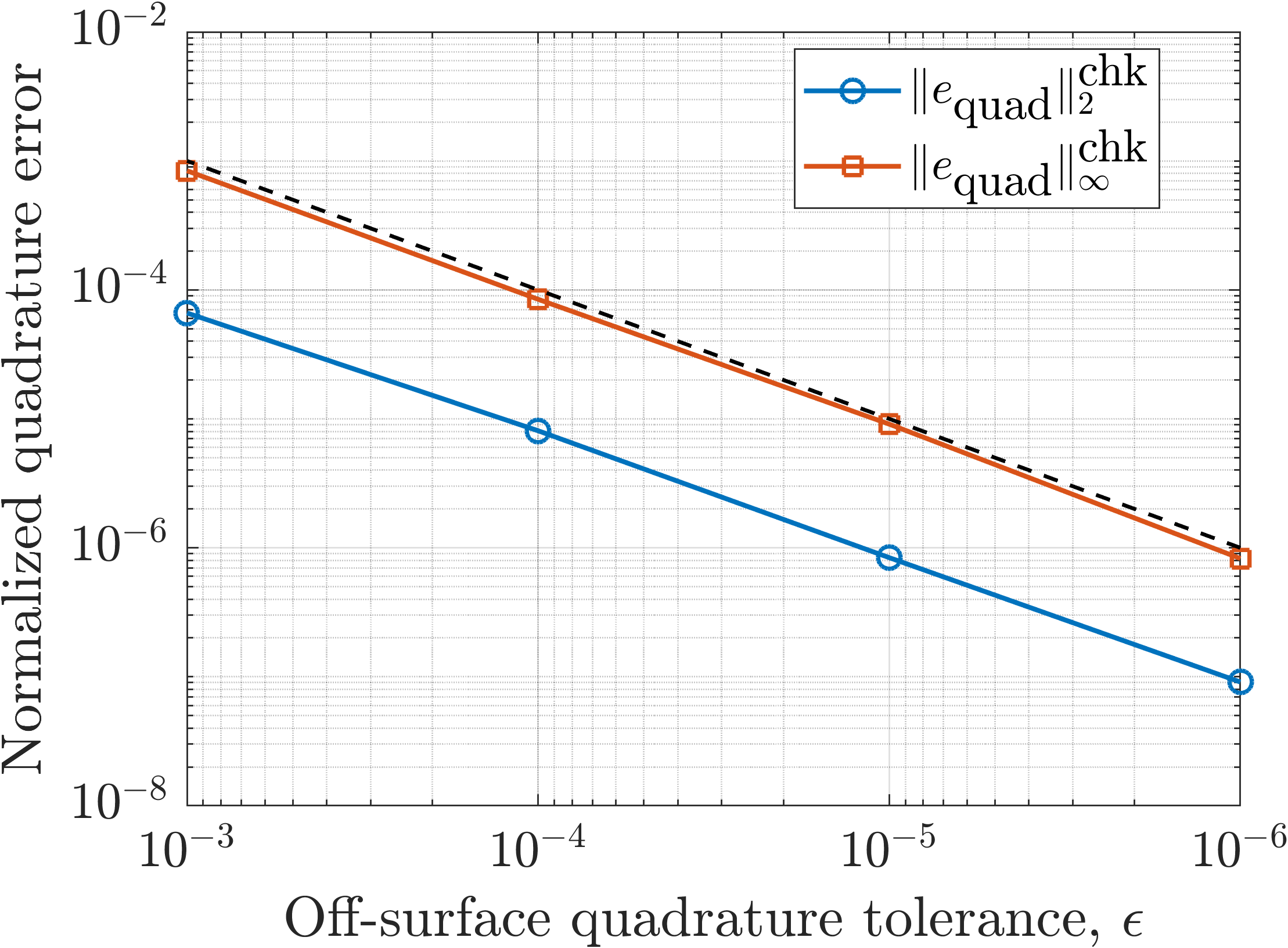}
\vspace{0.1em}
\caption{Quadrature error as a function of the off-surface quadrature tolerance $\epsilon$, evaluated over the five check surfaces shown in Figure \ref{fig:cluster_residual}. The dashed reference line shows where the quadrature error equals $\epsilon$.}
\label{fig:cluster_error_vs_tolerance}
\end{minipage}
\end{figure}

Figure \ref{fig:cluster_residual} shows $\log_{10}(e_{\mathrm{ns}})$ on the five check surfaces. The largest errors typically occur in the narrow gaps between neighboring particles. Over these surfaces, $\|e_{\mathrm{ns}}\|_{2}^{\mathrm{chk}}=2.15\times10^{-7}$ and $\|e_{\mathrm{ns}}\|_{\infty}^{\mathrm{chk}}=1.37\times10^{-5}$. Comparable error levels are obtained when other subsets of particles are checked, indicating that the five displayed particles are representative of the larger configuration. Thus, the no-slip condition is satisfied to about five digits in the maximum norm.

Next, we evaluate $\|e_{\mathrm{quad}}\|_{2}^{\mathrm{chk}}$ and $\|e_{\mathrm{quad}}\|_{\infty}^{\mathrm{chk}}$ for several values of $\epsilon$, while keeping the solved density, check points, DMK and QBX parameters, and available upsampling factors fixed. Figure \ref{fig:cluster_error_vs_tolerance} shows the results. The maximum error closely follows and remains just below the requested tolerance, while the area-normalized $L^2$ error is smaller still. 

\section{Conclusions}\label{s:conclusions}

Our two complementary goals have been to achieve a prescribed accuracy in the off-surface evaluation of layer potentials and to do so efficiently, without incurring unnecessary computational cost. To this end, we have developed a tolerance-driven quadrature framework for the Stokes double layer potential near smooth axisymmetric particles.

For each target--particle interaction, a fast target classifier selects the least costly quadrature treatment estimated to meet the prescribed tolerance: standard quadrature, the smallest adequate upsampling factor from a prescribed set, or special quadrature. Axisymmetry allows the geometry-dependent unit-density error indicators used by the classifier to be precomputed and tabulated in reduced cylindrical coordinates, while the local dependence on the layer density is incorporated through an inexpensive modifier on-the-fly. The resulting classification cost is negligible compared with that of the subsequent potential evaluation.

For interactions requiring special quadrature, we use a stabilized version of S3Q. The azimuthal integral is evaluated using TSSQ, which prevents the severe cancellation otherwise encountered for the Stokes double layer potential at very close targets. The associated target-specific azimuthal quadrature weights are constructed efficiently by one-dimensional interpolation from precomputed tables. In the polar direction, error indicators guide adaptive panel refinement of the Gauss--Legendre quadrature. Although the classifier can be coupled to other special quadrature methods, S3Q is particularly well suited to the present workflow because its automatic parameter selection allows the error to be controlled to the prescribed tolerance.

The off-surface framework has been integrated into a boundary integral solver for Stokes mobility and resistance problems. Precomputed QBX treats on-surface self-interactions \cite{AFKLINTEBERG2016}, while the classifier determines the quadrature treatment for off-surface target--particle interactions during both the boundary integral solve and subsequent flow evaluation. The numerical experiments range from detailed two-particle configurations to a triply periodic resistance problem involving $100$ closely spaced spheroids, in which the standard-quadrature contributions are accelerated using the periodic DMK method \cite{krantz2026dmk}. Across these experiments, the prescribed tolerance is met for nearly all tested target--particle interactions, and the error remains within a modest factor of the tolerance in the few remaining cases. 
The two-particle results further show that the classifier almost always selects the least costly quadrature treatment capable of attaining this accuracy. For example, at a tolerance of $10^{-6}$, only $0.02\%$ of the $59\,696$ interactions classified as requiring upsampling or special quadrature exceed the tolerance. By contrast, assigning these interactions to the next less expensive quadrature rule causes more than $95\%$ of them to exceed it; see \Cref{tab:s3q_errors,tab:cate_compare}. 

The solver used in this paper relies on a QBX implementation for on-surface evaluation that supports spheroids with a global tensor-product discretization. A natural direction for future work is to extend the complete solver to other smooth axisymmetric geometries, particularly those for which a panel-based polar discretization is advantageous. Both the fast target classifier and S3Q already extend to such settings, with the latter demonstrated in \cite{krantz2026s3q}. Although the present work focuses on the Stokes double layer potential, the off-surface framework can also easily be adapted to other Stokes layer potentials.

\bmhead{Acknowledgements}Matharu acknowledges support from the Verg foundation. Krantz and Tornberg acknowledge support from the Swedish Research Council under grant 2023-04269.

\begin{appendices}

\section{Numerical evaluation of one-dimensional error indicator integrals}\label{a:error_formulas}

Here we provide the formulas used in our implementation to evaluate the one-dimensional integrals appearing in the surface error indicators. These formulas were derived in \cite[Section~6.1]{SORGENTONE2023}, building on the earlier treatment in \cite[Section~6.1]{AFKLINTEBERG2022}.

For a general singularity order $p$, the directional error indicators in \eqref{eq:surface_unit_indicator_phi} and \eqref{eq:surface_unit_indicator_theta} contain the integrals
\begin{equation}
\begin{split}
J_{\contphi}(\xtil)
&\coloneqq
\int_0^\pi
\Psi_{n_\phi,p}^{\mathrm{TZ}}
\bigl(\contphi_0(\conttheta)\bigr)
\dif\conttheta,\\
J_{\conttheta}(\xtil)
&\coloneqq
\int_0^{2\pi}
\Psi_{n_\theta,p}^{\mathrm{GL}}
\bigl(t_0(\contphi)\bigr)
\dif\contphi,
\end{split}
\label{eq:Jintegrals}
\end{equation}
where $\Psi_{n_\phi,p}^{\mathrm{TZ}}$ and $\Psi_{n_\theta,p}^{\mathrm{GL}}$ are the one-dimensional quadrature factors defined in \eqref{eq:tz_root_factor} and \eqref{eq:gl_root_factor}, respectively.
The functions $\contphi_0(\conttheta)$ and $\conttheta_0(\contphi)$ are the selected roots of the squared-distance function in the azimuthal and polar directions, respectively, and $t_0(\contphi)=(2/\pi)\conttheta_0(\contphi)-1$ is the polar root mapped to the Gauss--Legendre reference interval $[-1,1]$. 

The integrals in \eqref{eq:Jintegrals} are generally not available in closed form. We therefore approximate them by combining a locally corrected root approximation with a low-order Gauss--Laguerre quadrature.

Let $(\bar\theta,\bar\phi)$ be the surface node closest to the target $\xtil$, and define
\begin{equation*}
\bar\x=\bgamma(\bar\theta,\bar\phi),\qquad 
\vec{r}=\bar\x-\xtil,\qquad
\bar\bgamma_\theta=\partial_\conttheta\bgamma(\bar\theta,\bar\phi),\qquad
\bar\bgamma_\phi=\partial_\contphi\bgamma(\bar\theta,\bar\phi),
\end{equation*}
Let $\bar\theta_0=\conttheta_0(\bar\phi,\xtil)$ be the corresponding root in the $\vartheta$-direction. The local approximation of the root in the $\conttheta$-direction is obtained from the local quadratic approximation of the squared-distance function in \eqref{eq:R2}, giving
\begin{equation*}
\theta_0^L(\contphi) = \bar\theta - \frac{b(\Delta\contphi)}{2c} \pm i \frac{\sqrt{4a(\Delta\contphi)c-b(\Delta\contphi)^2}}{2c},\qquad \Delta\contphi=\contphi-\bar\phi,
\end{equation*}
where
\begin{equation*}
\begin{split}
a(\Delta\contphi) &= \|\vec{r}\|^2 + 2(\vec{r}\cdot\bar\bgamma_\phi)\Delta\contphi + \|\bar\bgamma_\phi\|^2\Delta\contphi^2, \\
b(\Delta\contphi) &= 2(\vec{r}\cdot\bar{\bgamma}_\theta) + 2(\bar{\bgamma}_\phi\cdot\bar\bgamma_\theta)\Delta\contphi,\\
c                 &= \|\bar\bgamma_\theta\|^2.
\end{split}
\end{equation*}
To better capture the peak magnitude of the error indicator, this local approximation is corrected by the root value at the closest node. Thus, we define
\begin{equation*}
\widetilde\theta_0(\contphi) = \bar\theta_0 - \theta_0^L(\bar\phi) + \theta_0^L(\contphi),
\end{equation*}
where $\widetilde{t}_0(\contphi)=(2/\pi)\widetilde{\theta}_0(\contphi)-1$ is the corresponding mapped root. 
The same construction is used in the $\contphi$-direction, producing a corrected local approximation $\widetilde\phi_0(\conttheta)$.

Writing $K=\|\bar\bgamma_\theta\|/\|\bar\bgamma_\phi\|$, we use these corrected local roots to capture the exponential decay of the one-dimensional indicators away from the closest node. This gives the practical formulas
\begin{align*}
J_\contphi(\xtil)   &\approx \frac{1}{n_\phi K}\left(\int_0^\infty h_{-}^{\mathrm{TZ}}(s)e^{-s}\dif s+\int_0^\infty h_{+}^{\mathrm{TZ}}(s)e^{-s}\dif s\right),\\
J_\conttheta(\xtil) &\approx \frac{\pi K}{4n_\theta}\left(\int_0^\infty h_{-}^{\mathrm{GL}}(s)e^{-s}\dif s+\int_0^\infty h_{+}^{\mathrm{GL}}(s)e^{-s}\dif s\right),
\end{align*}
where
\begin{align*}
h_{\pm}^{\mathrm{TZ}}(s) &= \Psi_{n_\phi,p}^{\mathrm{TZ}}\left(\widetilde{\phi}_0\left(\bar\theta\pm\frac{s}{n_\phi K}\right)\right)e^s,\\
h_{\pm}^{\mathrm{GL}}(s) &= \Psi_{n_\theta,p}^{\mathrm{GL}}\left(\widetilde{t}_0\left(\bar\phi\pm\frac{\pi Ks}{4n_\theta}\right)\right)e^s.
\end{align*}
Each of the four half-line integrals is then evaluated using an 8-point Gauss--Laguerre rule. Since the indicators are used only for target classification, this low-order auxiliary quadrature is sufficient in practice.

\section{Analytical formulas for roots of the squared-distance function}\label{a:roots}
Here we collect the analytical formulas for the complex-valued roots of the squared-distance function used in the error indicators from Section \ref{s:precomputed_error_indicators}. We consider an axisymmetric surface $\Gamma$ parametrized by $\bgamma(\conttheta,\contphi)$ as in \eqref{eq:axisymmetric_param},
\begin{equation*}
\boldsymbol{\gamma}(\conttheta,\contphi) = \big(a(\conttheta) \, \sin(\conttheta) \cos(\contphi),~a(\conttheta) \, \sin(\conttheta) \sin(\contphi),~c(\conttheta) \, \cos(\conttheta)\big),
\end{equation*}
and a target point $\xtil=(x,y,z)\notin\Gamma$, with cylindrical radius $\rho=\sqrt{x^2+y^2}>0$. The associated squared-distance function is $R_{\bgamma}^2(\conttheta,\contphi,\xtil)$, defined in \eqref{eq:R2}. For a general axisymmetric surface, the root in the $\contphi$-direction is available in closed form, whereas the root in the $\conttheta$-direction is available analytically only in special cases such as the sphere and the spheroid.

\begin{lemma}[Root of $R_{\bgamma}^2$ in $\contphi$ for fixed $\conttheta$, general axisymmetric surface {\cite[Lemma~2]{SORGENTONE2023}}]
For a fixed $\conttheta\in(0,\pi)$, the equation $R_{\bgamma}^2(\conttheta,\contphi,\xtil)=0$ has the complex roots
\begin{equation*}
\contphi_0(\conttheta) = \atantwo(y,x) \pm i\ln\left(\frac{\lambda(\conttheta)+\sqrt{\lambda(\conttheta)^2-\rho^2\sin^2(\conttheta)}}{\rho\sin(\conttheta)}\right),
\label{eq:phi0}
\end{equation*}
where
\begin{equation*}
\lambda(\conttheta) = \frac{\tilde{a}(\conttheta)^2+\rho^2+\big(\tilde{c}(\conttheta)-z\big)^2}{2a(\conttheta)}, \quad\tilde{a}(\conttheta)=a(\conttheta)\sin(\conttheta),\quad\tilde{c}(\conttheta)=c(\conttheta)\cos(\conttheta),
\label{eq:lambda}
\end{equation*}
with $\lambda(\conttheta)>\rho\sin(\conttheta)$. Here, and in the following, $\atantwo(\eta,\xi)$ denotes the argument of the complex number $\xi+i\eta$, $-\pi<\atantwo(\eta,\xi)\leq\pi$. We also define $\sqrt{\zeta^2-1}$ as $\sqrt{\zeta+1}\sqrt{\zeta-1}$ with branch cuts chosen so that $-\pi<\arg(\zeta\pm1)\leq \pi$.
\end{lemma}

\begin{lemma}[Root of $R_{\bgamma}^2$ in $\conttheta$ for fixed $\contphi$, sphere {\cite[Lemma~3]{SORGENTONE2023}}]
Let $a(\conttheta)=c(\conttheta)=a$, and fix $\contphi\in[0,2\pi)$. Assume that if $z=0$, then $\contphi-\atantwo(y,x)\neq \pi/2+p\pi$, $p\in\mathbb{Z}$. Then the equation $R_{\bgamma}^2(\conttheta,\contphi,\xtil)=0$ has the complex roots
\begin{equation*}
\conttheta_0(\contphi)
=
\atantwo(x\cos\contphi+y\sin\contphi,\;z)
\pm
i\ln\!\left(
\lambda_0(\contphi)+\sqrt{\lambda_0(\contphi)^2-1}
\right),
\label{eq:theta0_sphere}
\end{equation*}
where
\begin{equation*}
\lambda_0(\contphi)
=
\frac{a^2+\rho^2+z^2}
{2a\sqrt{(x\cos\contphi+y\sin\contphi)^2+z^2}},
\end{equation*}
with $\lambda_0(\contphi)>1$.
\end{lemma}

\begin{lemma}[Root of $R_{\bgamma}^2$ in $\conttheta$ for fixed $\contphi$, spheroid {\cite[Lemma~6.2]{krantz2026s3q}}]
Let $a(\conttheta)=a$, $c(\conttheta)=c$ with $a\neq c$, and fix $\contphi\in[0,2\pi)$. Assume that if $z=0$, then $\contphi-\atantwo(y,x)\neq \pi/2+p\pi$, $p\in\mathbb{Z}$. Then the equation $R_{\bgamma}^2(\conttheta,\contphi,\xtil)=0$ has the complex roots
\begin{equation*}
\conttheta_0(\contphi)=\Arg(\beta(\contphi))-i\ln|\beta(\contphi)|,
\label{eq:theta0_spheroid}
\end{equation*}
where $\beta$ satisfies the quartic equation
\begin{equation*}
\frac{\Delta}{4}\beta^4+\tau\beta^3+\left(\frac{\Delta}{2}+d^2\right)\beta^2+\tau^*\beta+\frac{\Delta}{4}=0,
\end{equation*}
with
\begin{equation*}
\Delta=c^2-a^2,
\qquad
\tau=-cz+ia\left(x\cos\contphi+y\sin\contphi\right),
\qquad
d^2=a^2+\rho^2+z^2.
\end{equation*}
Among the four roots of this quartic, the relevant one is the one yielding the smallest $|\Im(\conttheta_0)|$.
\end{lemma}

\end{appendices}

\begingroup
\setlength{\bibsep}{8pt}
\bibliography{references}
\endgroup

\end{document}